\documentclass{article}
\usepackage{url}

\usepackage{setspace}  
\usepackage{pdfscreen}
\usepackage{hyperref}
\usepackage{mathptmx}
\usepackage{latexsym}
\usepackage{amssymb}
\usepackage{amsthm}
\usepackage{amsmath}
\usepackage{makecell}

\newtheorem{theorem}{Theorem}
\newtheorem{lemma}{Lemma}
\newtheorem{corollary}{Corollary}

\newtheorem{example}{Example}

\begin{document}                                                                   
\title{\huge Irreducible triangulations of the once-punctured torus}

\author{\Large S. Lawrencenko, T. Sulanke,  M.\,T. Villar, \\  \Large L.\,V. Zgonnik, M.\,J. Ch\'avez, J.\,R. Portillo}
\date{}
\maketitle

\begin{abstract}
\noindent A triangulation of a surface with fixed topological type is called irreducible if no edge can be contracted to a vertex while remaining in the category of simplicial complexes and preserving the topology of the surface. A complete list of combinatorial structures of irreducible triangulations is made by hand for the once-punctured torus, consisting of exactly 297 non-isomorphic triangulations.
\end{abstract}



\section{Introduction}
\label{S1}

Let $S \in \{ S_h , N_k \}$
be a {\it closed surface}---that is, the closed orientable (connected compact) $2$-manifold  $S_h $   of genus   $h$
  or the closed nonorientable $2$-manifold    $ N_k$  of nonorientable genus $k$.  Using this terminology, $S_0$  is the
  sphere, $S_1 $ is the torus, $N_1$  is the
 projective plane and $N_2$ is the Klein bottle. Let  $D$ be an open disk in  $S$, with boundary
 $\partial D = \partial (S-D) $    homeomorphic to a circle. In
 particular,  $S_0 -D$ is a disk, $N_1 -D$  is the M\"{o}bius band, and   $S_1 -D$  is the punctured torus. We use
 the notation $``\Sigma$''  whenever we assume the general case in which
 $\Sigma$ is meant to be either $S$ or $S-D$.

 If a graph $G$  is $2$-cell embedded in  $\Sigma$, the components of
  $\Sigma -G$ are called {\it faces}. A {\it triangulation} of $\Sigma$
  with a simple graph $G$  (where ``simple'' means  ``without loops and without parallel edges'') is a $2$-cell
 embedding  $T: G \to \Sigma$  in which each face is bounded by a {\it 3-cycle},
  that is, a cycle of length $3$ made up of $3$ vertices connected by $3$ edges of $G$; moreover,
  we demand that the closures of any two faces are either disjoint, share a single vertex, or share a single edge.  We denote by
  $V=V(T)$, $E=E(T)$, and  $F=F(T)$ the sets of vertices, edges, and faces
   of  $T$, respectively. The cardinality  $\mid V(T) \mid$ is called the
   {\it order} of $T$. By $G(T)$  we denote the graph $(V(T),E(T))$  of  $T$. Two triangulations  $T_1$ and $T_2$ are called
   {\it isomorphic} (denoted  $T_1 \cong T_2$) when there exists a bijection, called an {\it isomorphism},
   $\varphi : V(T_1) \to V(T_2)$, such that  $[u,v,w] \in F(T_1)$ if and only if
   $[\varphi(u), \varphi (v), \varphi (w)] \in F(T_2)$. Throughout this paper we
   distinguish  between triangulations only up to isomorphism. For $\Sigma =S-D$,
    let  $\partial T$ ($= \partial D$) denote the boundary cycle of  $T$.
    The vertices and edges of $\partial T$  are called {\it boundary 
    vertices} and {\it boundary edges} of  $T$.

A triangulation of a 2-manifold with fixed topological type is viewed as a member of the category
of simplicial complexes. A triangulation is called {\it irreducible} if no edge can be contracted
(to a vertex) without vacating the category of simplicial complexes or changing the topology of the underlying 
$2$-manifold. Obstacles for edge contraction are studied in Section \ref{S3}; one typical obstacle
is the creation  of parallel edges (forbidden in a simplicial complex).  The term ``irreducible
triangulation'' is more accurately introduced in Section \ref{S3}.  For the sake of brevity we
abbreviate ``irreducible triangulation'' as ``IT'', and ``irreducible triangulation of the (once-)
punctured torus'' as ``ITPT''.

\medskip
\medskip
    \begin{center}
    \hspace*{1cm}    \pdfimage width .85\textwidth {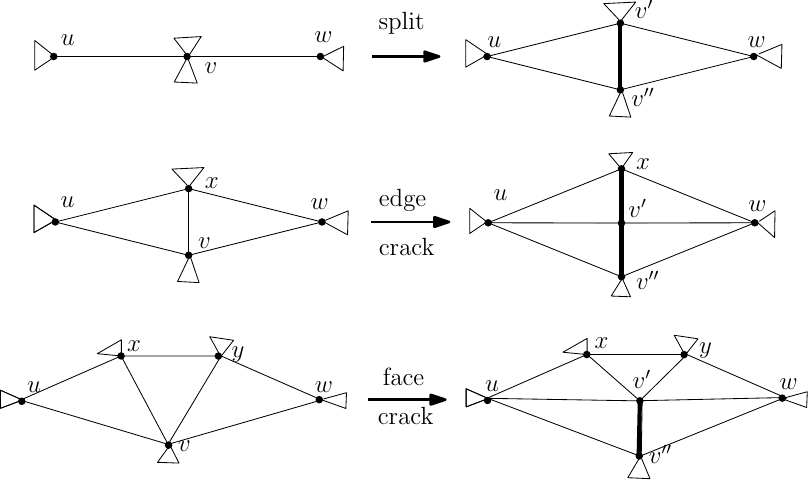}

    Figure 1. Splitting (top), edge cracking (middle), face cracking (bottom).
    \end{center}
\medskip
\medskip

The collection of all ITs of $\Sigma$   is a basis for the family of all triangulations of  $\Sigma$,  in the sense that any triangulation of  $\Sigma$ can be
obtained from a member of the basis by repeating the
{\it splitting} operation introduced in Section \ref{S3} and illustrated by Figure 1 (top) along with two special cases in the middle and bottom. Barnette and Edelson \cite{BE} and independently Negami
\cite{N} proved that for every closed $2$-manifold  $S$ the basis of ITs is finite.

     At present, bases of ITs are known for seven closed $2$-manifolds: the sphere (Steinitz \cite{SR}, Bowen and Fisk \cite{BF}),
      projective plane (Barnette \cite{B}), torus (Lawrencenko \cite{L, L-2}), Klein bottle (Lawrencenko and Negami \cite{LN},
       and Sulanke \cite{S}) as well as $S_2$, $N_3$, and   $N_4$ (Sulanke \cite{S-2, S-3, S-4}). In Section  \ref{S2}, we briefly
       consider the case of the sphere in the historical retrospect.

     In 2012 the first author (Lawrencenko) proposed \cite{L-5} the problem of determining all ITs of a given $2$-manifold with boundary,
     that is, all irreducible ``punctured'' triangulations with given genus and number of punctures (the latter is equal to the number
     of boundary components). This problem is motivated by the enormous growth of the number of ITs on closed $2$-manifolds as the genus
grows: already $396,784$ ITs on  $S_2$, and $6,297,982$ ITs on $N_4$  (see \cite[p. 3]{S-3}).
Firstly, Boulch, Colin de Verdi\`{e}re, and Nakamoto \cite{BCN} gave upper bounds on the order of
an IT of the punctured $2$-manifold with given genus and number of punctures; their bounds are
asymptotically tight but loose for low genus. Then, three of the authors of the present paper with
participation of Quintero \cite{CLQV, CLQV-2} produced a complete list of the $6$ ITs on
the M\"{o}bius band.

     Before the current study, a nearly comprehensive list of $293$ non-isomorphic ITPTs was already obtained with the aid of a computer program in \cite{CLPV} and was announced in \cite{L-8}. That list missed $4$ ITPTs and was corrected in \cite{CLPV-2}. In this paper we develop a fully comprehensive list of $297$ ITPTs (Theorem \ref{12.1}) independently and without using a computer. (The terminology used in \cite{CLPV, CLPV-2}   is slightly different from that used in this paper.) Throughout this paper we intentionally avoid the use of a computer and all work can be checked tediously by hand.

\section{The sphere}
\label{S2}

In this section we examine more closely the case of the sphere; it is a very
instructive one. It is not very hard to prove the following theorem; this is the content of Exercise 1  \cite[p. 243]{G}  and Exercise 6 \cite[p. 59]{EH}.

\begin{theorem}{\em(\cite{SR, BF})}.
\label{2.1}
There is only one irreducible triangulation of the sphere: the boundary of a tetrahedron.
\end{theorem}

In the historical setting, Bowen and Fisk \cite{BF} were the first who brought a modern version of
Theorem \ref{2.1}; in fact, \cite{BF} contains a stronger result. However, Theorem  \ref{2.1} can
be derived from the considerations in the book \cite[pp. 227-229]{SR} of Steinitz and Rademacher
(1934), although the authors use a different terminology. In that book it is shown that (after
appropriate translation and interpretation) every $3$-regular polyhedron can be obtained from a
tetrahedron by a succession of the face-splitting operations. This means in dual form that the
tetrahedral triangulation is a unique spherical IT.

The shortest proof  \cite{L-7} of Theorem \ref{2.1} uses a corollary of \cite[Theorem 2, p. 264]{LNS}. In fact, that corollary is a characteristic property of the sphere; it distinguishes the sphere from the rest of the $2$-manifolds. The corollary states that every triangulation of the sphere contains a {\it clean} vertex---that is, a vertex  $v$ whose link is chordless; the {\it link} of  $v$ consists of the cycle of vertices and edges surrounding  $v$. Then $v$ is not incident with any {\it non-facial} 3-cycle---that is, a 3-cycle that does not bound a face of the triangulation. Thus any edge incident with $v$  can be contracted without creating parallel edges. One can keep iterating the edge contraction process until it terminates at some triangulation with only four faces. Although all vertices are still clean in such a tetrahedral triangulation, no edge is contractible anymore, but for another reason: after having attempted to contract an edge in the tetrahedral triangulation, one gets to a doubly covered triangle, which is not a simplicial complex.

\section{Preliminaries}
\label{S3}
Let $T$  be a triangulation of $\Sigma$. An unordered pair of distinct
adjacent edges $[v,u], [v,w]$ of  $T$ is called a {\it divider} of  $T$  ({\it centered }) {\it at
vertex} $v$, denoted by $\langle u, v, w \rangle$ ($= \langle w, v, u \rangle$). The {\it splitting} of  $ \langle u, v, w \rangle$, denoted  ${\rm{sp}} \langle u,v,w \rangle$, 
is the operation which consists of deleting $\langle u, v, w \rangle$ from $T$ and
closing the resulting hole with two new triangular faces, $[v',v'',u]$ and $[v',v'',w]$, where  $v'$ and $v''$  denote the two split vertices;  see Figure 1 (top). Under this operation, the vertex $v$  is
 extended to the edge $[v,'v'']$  and the two faces incident with this edge are inserted into the triangulation. Specifically in the case in which
 ($\Sigma = S-D$) {\it AND} ($[u,v] \in E(T) $) {\it AND} ($ v \in V(\partial T) $) {\it AND} ($ u \notin V(\partial T) $), the
 operation sp$\langle u,v]$  of {\it splitting a truncated divider} $\langle u,v]$
   produces a single new triangular face $[u,v',v'']$, where  $[v',v''] \in E\left(\begin{array}{c}\partial \left( {\rm{sp}} \langle u,v](T)\right)
   \end{array}  \right)$.

Under the inverse operation, which is called {\it  contracting} the edge $[v',v'']$,  this edge collapses into a single vertex  $v$, the faces  $[v',v'',u]$  and  $[v',v'',w]$  collapse into single edges $[v,u]$
and  $[v,w]$, respectively. This operation is denoted by
 ${\rm{sh}} \rangle v', v'' \langle $, which comes from the word ``shrinking'', a synonym for ``contracting''.
 Therefore, ${\rm{sh}} \rangle v',v'' \langle \left( \begin{array}{c} {\rm{sp}}\left\langle  u, v, w\right\rangle (T) \end{array}\right) =T$. It
should be noticed that in the case
($ \Sigma=S-D)$ {\it AND} ($[v',v'']  \in E(\partial T) $), there is only one face
incident with $[v',v'']$ and that face collapses to a single edge
under  sh$\rangle v',v'' \langle$. Clearly, the operation of splitting
never changes the topology  of  $\Sigma$. We demand that the contraction operation
must preserve the topology of $\Sigma$  as well; moreover, parallel
 edges must not be created in a triangulation.  In the case in which
 an edge $\varepsilon \in E(T)$  occurs in some non-facial $3$-cycle, if we still
 insist on contracting  $\varepsilon$, parallel edges would be produced, which would
 exclude  ${\rm{sh}} \rangle \varepsilon \langle (T)$ from the category of triangulations. An edge $\varepsilon $
  is called {\it contractible} or a {\it rope} (or a{ \it cable}) if ${\rm{sh}} \rangle \varepsilon \langle (T)$
  is still a triangulation of  $\Sigma$; otherwise $\varepsilon$ is called
  {\it noncontractible} or a {\it rod}. Therefore, one can contract ropes but not rods. The subgraph of $G(T)$ made up of all ropes is called the {\it rope subgraph} of  $G(T)$.
  
The only constraint to edge-contractibility in a  non-tetrahedral triangulation $T$  of a
\textit{closed} $2$-manifold  $S$ is determined in
\cite{B, BE, L, L-2}: an edge  $\varepsilon \in E(T)$
is a rod if and only if $\varepsilon$ satisfies the following condition: \\

 (3.1) $\varepsilon$ is in a non-facial $3$-cycle of  $G(T)$. \\

That is, one cannot contract the edges of a non-facial 3-cycle.

     There are, in all, three constraints to edge-contractibility in a triangulation  $T$ of a {\it punctured}
$2$-manifold $S-D$. Two of them are determined in \cite{BCN}: an edge $\varepsilon \in E(T)$
  is a rod if and only if $\varepsilon$  satisfies either condition (3.1)
  or the following condition: \\

(3.2)  $\varepsilon$  is a {\it chord} of $D$---that is, the end vertices of  $\varepsilon$ are
in $V(\partial D)$  but  $\varepsilon \notin E(\partial D)$. \\

That is, one cannot contract chords.

A triangulation is said to be {\it irreducible} if it has no ropes, or
in other words, each edge is a rod. For instance, a single triangle is the
only IT of the disk $S_0-D$  although its edges don't
 meet either of conditions (3.1) and (3.2). Thus, the following is  yet one more
constraint to edge-contractibility of  $\varepsilon$: \\

(3.3) $\varepsilon$  is a boundary edge in the case the boundary cycle is a
$3$-cycle. \\

Although condition (3.3) is a specific case of condition (3.1)
 unless $S=S_0$  and is not explicitly stated in \cite{BCN}, this condition deserves special mention.
 In the remainder of this paper we assume that  $S \neq S_0$.

\section{The structure of irreducible punctured triangulations}
\label{S4}

 As an energy metaphor (with ropes thought of as high voltage cables), a vertex of a triangulation   $R$ of a closed $2$-manifold $S$ is called a {\it pylonic vertex} or a $\ast$-{\it  vertex}
  if that vertex is incident with all ropes of  $R$. A triangulation
  that has at least one rope and at least one $\ast$-vertex is called a
  {\it pylonic triangulation} or a $\ast$-{\it triangulation}. In other words, $R$  is pylonic if and only if the rope subgraph of  $R$ is isomorphic to the complete bipartite graph $K_{1,n}$  for some $n$.

Let $T$  be a triangulation of  $S-D$. Let us restore the disk  $D$ in  $T$,
 add a vertex $p$  in  $D$ and join $p$  to the vertices in $\partial D$.
We thus obtain a triangulation,  $T^*=T \cup D$, of $S$.
We call  $D$ the {\it patch}, call $p$ the {\it central vertex of the
patch}, and call $T^*$ a {\it parent triangulation}  for $T$.

 It will be shown in Section \ref{S10} that there exist single-roped triangulations of  $S_1$---that is, triangulations that have only one rope and thus two $\ast$-vertices. However, if a $\ast$-triangulation  $R$  has at least two ropes,  $R$ has a unique $\ast$-vertex. It is to be noted that if $T$  is an IT of $S-D$, then  $T^*$ may turn out to be an IT of  $S$, but not necessarily.

\begin{lemma}\label{4.1}
 If  $T$  is an IT of $S-D$, then unless $T^*$ is irreducible $T^*$ is pylonic.
\end{lemma}

\begin{proof}
Let $\varepsilon$ be an edge of $T$,  necessarily a rod. Then $\varepsilon$   satisfies either (3.1) or (3.2). If   $\varepsilon$  satisfies (3.1), the corresponding edge $\varepsilon^*$   in  $T^*$  also satisfies (3.1). If   $\varepsilon$ satisfies (3.2),  $\varepsilon^*$  still satisfies (3.1). In either case,  $\varepsilon^*$   is a rod. Thus, all ropes of   $T^*$ (if any) are not edges of  $T$ and therefore are incident with the central vertex of the patch.
\end{proof}

\begin{corollary}\label {4.2}
 A triangulation $T$ of $S-D$ is irreducible if and only if  $T$ is obtained from a parent triangulation  $T^*$  of $S$  either by deleting a vertex when  $T^*$ is irreducible, or by deleting a $\ast$-vertex   $p$ when  $T^*$ is pylonic.
\end{corollary}

\begin{proof}
The ``only if'' part follows from Lemma \ref{4.1}. The ``if'' part is trivial: any edge that is not
incident with  $p$ but occurs in some non-facial $3$-cycle through  $p$ is a chord of the link of
$p$ and therefore is still a rod after the deletion of $p$.
\end{proof}

     It is not a trivial question as to how the rope subgraph evolves under successively performed splittings. However,
      the following observation is easy to see: it is hard to make a rod out of a rope.

\begin{lemma}\label{4.3}
 The only situation in which a rope in a triangulation of $S$  changes into a rod
under a single application of the splitting operation is when the splitting is equivalent to the stellar subdivision
of a face having that rope as a boundary edge.
\end{lemma}

\begin{proof}
In such a situation, the edge contraction inverse to the splitting changes some rod $\varepsilon$
into a rope. Since the contraction does not change the homotopy type of the cycles in $\pi_1(S)$,
it follows by condition (3.1) above that  $\varepsilon$ occurs in a non-facial, null-homotopic
$3$-cycle, which is possible if and only if the disk bounded by that $3$-cycle is stellar
subdivided.
\end{proof}

\begin{corollary}\label {4.4}
If a triangulation of  $S$ is neither irreducible nor pylonic, it can  never become pylonic after any sequence of splittings,
except the case in which the rope subgraph is a $3$-cycle forming the boundary of a face.
\end{corollary}

We refer to the exceptional case of Corollary \ref{4.4} as a ``$\Delta$''. It is easy to see that a
$\Delta$  cannot occur in a once-split IT, but may arise in a twice-split IT of  $S$. By the way,
this is a good exercise for the reader to find an example of a twice-split IT  $K_{12} \to S_6$ in
which a $\Delta$  occurs. Interestingly, $\Delta$  does not occur for the Klein bottle but does
occur for  $S_2$ and  $N_3$, which is obtained by searching ITs with the second author's computer
program {\it surftri} \cite{S-4}. The example with $\Delta$  for  $S_2$ can be joined to an IT of  $S_1$
to produce an example with $\Delta$  for  $S_3$, etc. So $\Delta$ does occur for $S_h$  $(h \geq 2
)$ and, similarly, for $N_k$  $(k \geq 3 )$.

\begin{lemma}\label{4.5}
No $\Delta$ can ever be on the torus.
\end{lemma}

The proof of Lemma \ref{4.5} is postponed until the end of Section \ref{S11} when we will have more factual material to draw upon. By Lemma \ref{4.5}, we can restate Corollary \ref{4.4} as follows.

\begin{corollary}\label {4.6}
If a triangulation of  $S_1$ is neither irreducible nor pylonic, it can  never become pylonic after any sequence of splittings.
\end{corollary}

\section{The torus}
\label{S5}

Throughout the remainder of this paper, we only consider triangulations of $S_1$ or  $S_1-D$.

A theorem of the first author  \cite{L, L-2} states that
for $S_1$there exist, in all, twenty-one non-isomorphic ITs: 
${\rm{T}}_1, {\rm{T}}_2, \dots , {\rm{T}}_{21}$.
They are represented in each of Figures 2, 3, and 4 with their vertices numbered by $1, 2, \ldots,10$; each ${\rm{T}}_i$  is identical for the three figures, with fixed vertex numbering.  For each rectangle identify each pair of opposite sides to obtain an actual triangulation of  $S_1$.

An {\it automorphism} of a triangulation  $T$ is an isomorphism of $T$  with itself.
The set of dividers of $T$  as well as the sets $V(T), E(T), F(T)$   naturally fall into disjoint orbits under the action of the automorphism group  Aut($T$).
The  groups Aut({\rm{T}}$_i$) are determined explicitly for each    $i= 1, 2, \ldots, 21$ in  \cite{L-3}
 and are reproduced here in Table I in the form of generating sets. In particular, the generating set for
 Aut({\rm{T}}$_2$) was found in \cite[p. 544]{L-4}. Originally, the technique used in \cite{L-3} is based on a computer program,
  but it is a good exercise for the reader to verify without a computer that the results in Table I are correct.
   (Interestingly, the reader may notice from Figure 2 that for $i=6$  to  $17$, by appropriately
   dissecting the quadrilaterals into triangles, each  ${\rm{T}}_i$ can be obtained from the Cartesian product of two $3$-cycles quadrangularly embedded in  $S_1$; 
   the quadrangulation itself has a flag-transitive automorphism group of order $8 \times 3^2=72$; see \cite{L-6, L-9, LS}.)

Based on Table I, we have identified the vertex, face, and edge orbits in each ${\rm{T}}_i$  $(i=1,2, \dots,21)$.
Elements in the same orbit are called {\it similar}. Figure 2 shows the vertex orbits in each
 ${\rm{T}}_i$, where two vertices are marked by the same letter provided that the vertices are similar. Analogously,
 Figures 3 and 4 show the face and edge orbits, respectively.
 The same set of letters $\{ a,b,c,d, \dots \}$  is used for marking Figures~2, 3, and 4
 in which, for each $i$,
 the three sets $V({\rm{T}}_i)$, $E({\rm{T}}_i)$, and  $F({\rm{T}}_i)$, respectively, are marked up independently of each other. 

  \begin{center}
   \hspace*{0.0cm}     
\pdfimage width 1.0\textwidth {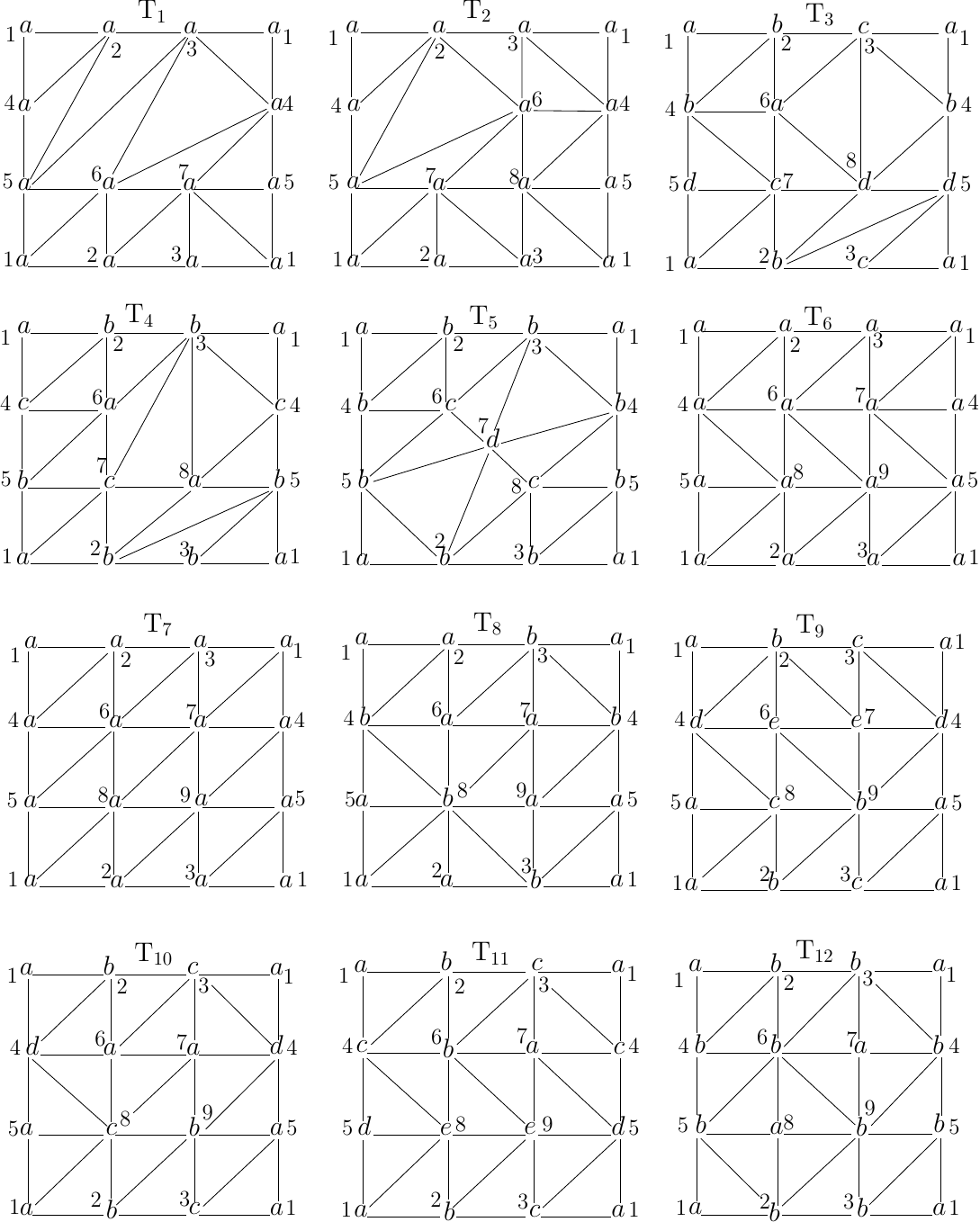}
\newline

     Figure 2. The vertex orbits (to be continued).
  \end{center}

\medskip

To count the number of vertex,
 face, or edge orbits in a specified  ${\rm{T}}_i$, we just count the number of distinct letters used for marking   ${\rm{T}}_i$
 in Figures 2, 3, or 4, respectively.

In what follows we implicitly use the obvious fact that if two dividers  $\langle u_1, v_1, w_1  \rangle$ and
 $\langle u_2, v_2, w_2  \rangle$   of ${\rm{T}}_i$  are similar, then the triangulations
     ${\rm{sp}}  \langle u_1, v_1, w_1  \rangle ({\rm{T}}_i)$ and ${\rm{sp}}  \langle u_2, v_2, w_2  \rangle ({\rm{T}}_i)$  are isomorphic.
     However, the converse is not generally true, as can be seen in the forthcoming sections; in particular,
     many counterexamples can be found in Table II.

\medskip

        \begin{center}
\hspace*{0.0cm}              
\pdfimage width 1.0\textwidth {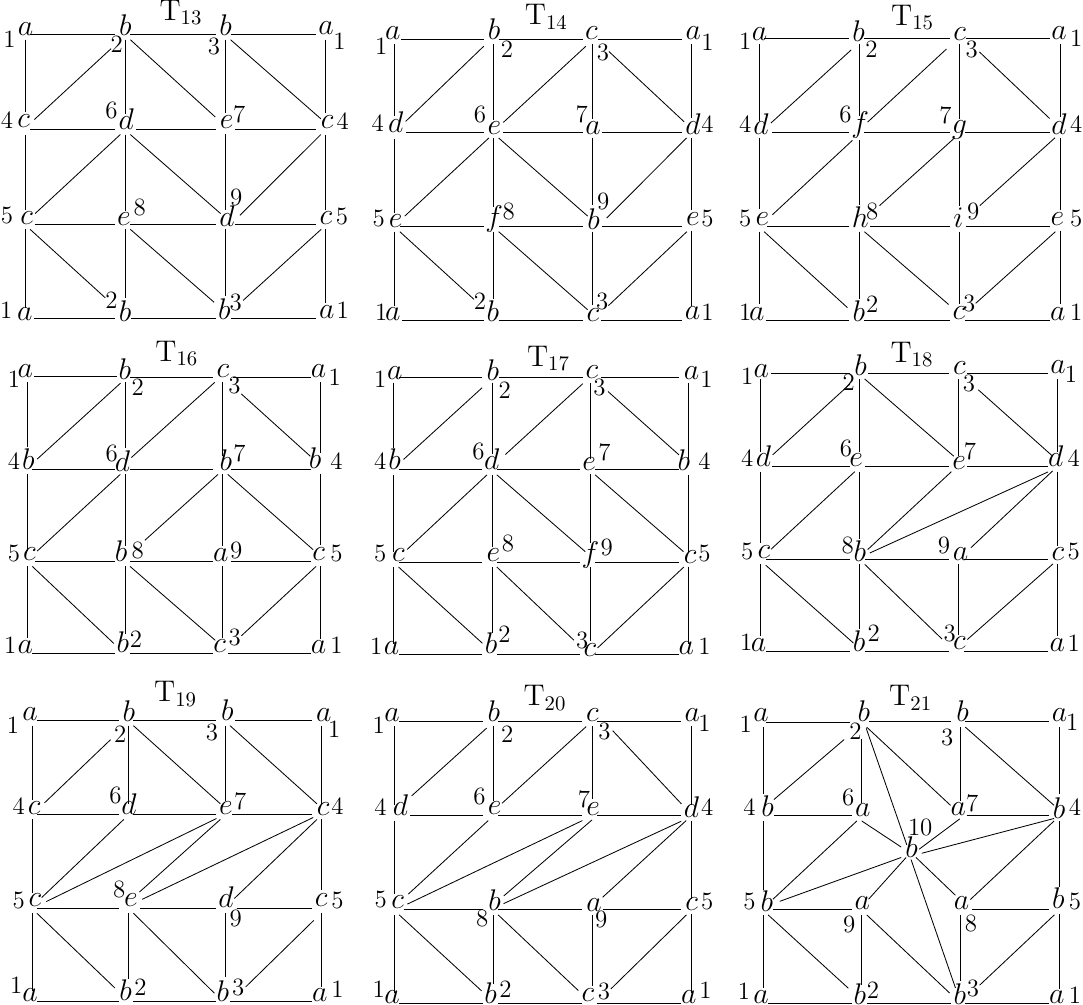}
\newline

        Figure 2. The vertex orbits (contd.)
     \end{center}

\medskip

\medskip

Let $T$  be an arbitrary triangulation of  $S_1$.  We define the {\it spread of a divider}  $\langle u, v, w  \rangle$  in  $T$,
 denoted $\mid u, v, w \mid$,  to be the {\it least} distance between $u$  and $w$  in the link of $v$. A divider with spread
  $k$  is referred to as a  $k$-divider. We observe from Figure 2 that the largest degree of a vertex in any
  ${\rm{T}}_i$ is $8$, thus the largest spread of a divider in any  ${\rm{T}}_i$ is $4$, and thus the set of all triangulations obtainable by single
  splitting from toroidal ITs can be written as $\Lambda= \bigcup_{k=1}^4 \Lambda_k,$ where
\begin{equation}
    \Lambda_k = \bigcup_{i=1}^{21} \{{\rm{sp}} \langle u, v, w  \rangle ({\rm{T}}_i) \, \big| \mid u,v,w \mid=k\}
\end{equation}

The spread of the divider used in generating the splitting is not an invariant of the triangulation obtained
because $\Lambda_3 \cap \Lambda_4\neq \emptyset$; for example, ${\rm{sp}}  \langle 2, 6, 9 \rangle ({\rm{T}}_6) \cong {\rm{sp}} \langle 1, 4, 5\rangle ({\rm{T}}_9)$ with  $\mid 2,6,9 \mid=3$ in  ${\rm{T}}_6$ and   $\mid 1,4,5\mid=4$ in  ${\rm{T}}_9$. The bottom line, however, is that we generate {\it all} triangulations of $S_1$  by repeatedly applying the splitting operation to the basis triangulations  ${\rm{T}}_i$.

\begin{lemma}\label{5.1}
$\Lambda_1 \cap \Lambda_2 = \emptyset$, $\Lambda_2 \cap \Lambda_3 = \emptyset$, and
 $\Lambda_3 \cap \Lambda_1 = \emptyset$.
\end{lemma}

\begin{proof}
For $k=1$  or $2$, the  $\ast$-vertex of any  $\ast$-triangulation $T$  in  $\Lambda_k$  is incident with at least two ropes and has degree  $3$ or  $4$, respectively, but any  $\ast$-triangulation in   $\Lambda_3$ has either only one rope, or at least two ropes with the degree of the  $\ast$-vertex of at least $5$, and the statement follows immediately.
\end{proof}

\begin{corollary}\label {5.2}
Any two ITPTs obtained by deleting a $\ast$-vertex from parent $\ast$-triangulations in different sets  $\Lambda_1$, $\Lambda_2$, or  $ \Lambda_3 $, are non-isomorphic.
\end{corollary}

\begin{proof}
We assume to the contrary that two such ITPTs $T$  and  $R$ are isomorphic. Then any isomorphism from $T$  to  $R$  takes  $\partial T$  onto  $\partial R$   and naturally extends to an isomorphism between the corresponding toroidal $\ast$-triangulations   $T^*$ and  $R^*$, which contradicts Lemma \ref{5.1}.
\end{proof}

\newpage
            \begin{center}
Table I. The automorphism groups of the toroidal ITs.

\scalebox{1.0}{
                \begin{tabular}{|c|c|c|}
                    \hline      &                                                  & \\
                    T$_i$&  Set of permutations generating Aut(T$_i$)         & $|$Aut(T$_i$)$|$   \\
                    &                                                  & \\
                    \hline T$_1$    &(1  5  3  4  7  2),  (1  2  3  4  5  6  7)        &    42  \\
                    \hline T$_2$    & (3 5)(4 7),  (1 6)(3 7)(4 5),  (1 5 2 7 6 3 8 4) &    32  \\
                    \hline T$_3$    & (2 4)(3 7),  (1 6)(3 7)(5 8)                     &    4   \\
                    \hline T$_4$    & (3 5)(4 7)(6 8),  (1 6 8)(5 3 2)                 &    6   \\
                    \hline T$_5$    & (2 3)(4 5)(6 8),  (2 5)(3 4)(6 8)                &    4   \\
                    \hline T$_6$    & (2 3)(4 5)(6 9)(7 8),  (1 6 5 2 7 8 3 4 9)       &    18  \\
                    \hline T$_7$    & (2 8 5 3 7 4)(6 9),  (4 8)(5 7)(6 9),            &        \\
                    & (1 2 3)(4 6 7)(5 8 9),  (1 5 4)(2 8 6)(3 9 7)    &    108 \\
                    \hline T$_8$    & (2 5)(3 4)(6 9), (1 5)(2 9)(3 8)(6 7),           &        \\
                    & (1 9 6)(2 5 7)(3 8 4)                            &    12  \\
                    \hline T$_9$    &   (1 5)(2 9)(3 8)(6 7)                           &    2   \\
                    \hline T$_{10}$   & (1 5)(2 9)(3 8)(6 7), (1 6)(3 8)(5 7)            &    4   \\
                    \hline T$_{11}$   & (1 7)(2 6)(3 4)(8 9)                             &    2   \\
                    \hline T$_{12}$   & (1 7 8)(2 4 9)(3 6 5),  (2 3)(4 5)(6 9)(7 8),    &        \\
                    & (2 4)(3 5)(7 8)                                  &    12  \\
                    \hline T$_{13}$   & (2 3)(4 5)(6 9)(7 8)                             &    2   \\
                    \hline T$_{14}$   & (1 7)(2 9)(5 6)                              &    2   \\
                    \hline T$_{15}$   &       The trivial automorphism group             &    1   \\
                    \hline T$_{16}$   & (2 4)(3 5)(7 8),  (1 9)(2 7)(4 8)                &    4   \\
                    \hline T$_{17}$   & (2 4)(3 5)(7 8)                                  &    2   \\
                    \hline T$_{18}$   & (1 9)(2 8)(3 5)(6 7)                             &    2   \\
                    \hline T$_{19}$   & (2 3)(4 5)(6 9)(7 8)                             &    2   \\
                    \hline T$_{20}$   & (1 9)(2 8)(3 5)(6 7)                             &    2   \\
                    \hline T$_{21}$   & (2 4 3 5)(6 7 8 9),  (1 6 7 9)(3 5 4 10),        &        \\
                    &             (1 8 6 9)(2 3 4 10)                  &    20  \\
                    \hline

                \end{tabular}
}
                
            \end{center}

 \medskip

 \medskip

Thanks to Corollary \ref{4.2}, the search for ITPTs is reduced to the search for vertex orbits in toroidal ITs
 and the search for toroidal $\ast$-triangulations, both in and out of  $\Lambda$. This task naturally splits into
 five cases depending on the origin of the parent triangulation $T$  in Corollary \ref{4.2}.

\section{The search for ITPTs: Case 0: Series 1}
\label{S6}

{\it Case 0. Parent triangulation} $T^*$ {\it   is an irreducible triangulation of} $S_1$. \\

     With help from Figure 2, we can count a total of  $80$ vertex orbits in  $T^* = {\rm{T}}_i$, where  $i$ runs over the set  $\{1, \dots, 21\}$. By deleting an arbitrary vertex in each of the orbits, we obtain Series 1 of $80$ non-isomorphic ITPTs, thanks to Corollary \ref{4.2}. This Series contains no isomorphic pairs as can be proved using an argument similar to the one used in Corollary \ref{5.2} by reduction to a contradiction due to the non-similarity of the vertices deleted.

  \section{The search for ITPTs: Case 1: Series 2}
\label{S7}

{\it Case 1. Parent triangulation} $T^*$ {\it   is in} $\Lambda_1$ {\it or can be obtained from a member of }  $\Lambda_1$
{\it by a sequence of splittings.} 
\medskip

     By Eq. (1), if $T^* = {\rm{sp}}   \langle  u, v, w \rangle  ({\rm{T}}_i) \in \Lambda_1$, then  $   \langle  u, v, w  \rangle $ is a $1$-divider and  ${\rm{sp}}  \langle  u, v, w \rangle  $
     is equivalent to the stellar subdivision of the face  $[ u, v, w]$. Thus, such a triangulation $T^*$ is pylonic with the only
     (necessarily $3$-valent) $\ast$-vertex---either  $v'$ or  $v''$ as a matter of notation. It can be easily seen that there are not
     any other $\ast$-triangulations that belong to  $\Lambda_1$  and that any further splitting of  $T^*$ would lead to a triangulation
     that is no longer pylonic. By Corollary \ref{4.6}, there are not any other $\ast$-triangulations obtainable from a member of
     $\Lambda_1$   by any sequence of splittings.

     The deletion of the 3-valent $\ast$-vertex from  $T^*$  is equivalent to the deletion of the corresponding face from
      ${\rm{T}}_i$. Figure 3 shows the face orbits---there are totally $129$ non-similar faces in  ${\rm{T}}_i$
      $( i= 1, \dots, 21)$. By deleting an arbitrary face in each of the $129$ orbits from  ${\rm{T}}_i$, we obtain Series~2 of 129
      non-isomorphic ITPTs thanks to Corollary \ref{4.2}; this Series is complete by Corollary \ref{4.6}.

     Just as in Section \ref{S6} it can be shown that the set of the $129$ triangulations contains no isomorphic pairs.

 \section{The search for ITPTs: Case 2: Series 3}
\label{S8}

{\it Case 2.}   Parent triangulation  $T^*$  is in $ \Lambda_2$ or can be obtained from a member of   
$\Lambda_2$ by a sequence of splittings.

     By Eq. (1), if  $T^* = {\rm{sp}}  \langle  u, v, w  \rangle  ({\rm{T}}_i) \in \Lambda_2$, then   
$\langle  u, v, w \rangle $ is a $2$-divider.
     Denote by $x$  the vertex determining a path of length $2$ in the link of  $v$ together with vertices
      $u$ and  $w$; see middle left of Figure 1. This specific type of splitting is equivalent to the
       {\it cracking of the edge } $[v,x]$---that is, adding a vertex  $v'$ to  $[v,x]$ and connecting $v'$  to the apices
         $u$  and $w$  of the triangular faces incident with $[v,x]$.  This transformation always leads to a triangulation
         with a new $4$-valent vertex  $v'$  at which the two ropes $[v',v'']$  and  $[v',x]$  form a 2-divider shown in bold in the
          middle right of Figure 1. Sometimes $v'$  may turn out to be pylonic, in which event it can be easily seen
          that any further splitting of   $T^*$ would lead to a triangulation that is no longer pylonic. Thus, by Corollaries \ref{4.2}
          and \ref{4.6}, our first goal is to find all $\ast$-triangulations   $T^*$ in  $ \Lambda_2$. 
Figure 4 shows the edge orbits in  ${\rm{T}}_i$. There are totally $203$ non-similar edges, but only $89$ of them, when cracked, actually produce $\ast$-triangulations, as we have checked by direct inspection; the $89$ edges can be seen in Table II.

The removal of the $4$-valent $\ast$-vertex is equivalent to the removal of the corresponding edge
(the one being cracked) from  ${\rm{T}}_i$ which produces a quadrilateral hole. It is easy to see that the removal of
similar edges gives isomorphic triangulations (of $S_1-D$), but it may happen that the removal of non-similar edges
produces isomorphic triangulations. It can be verified straightforwardly that there are exactly $27$ isomorphic pairs
among the $89$ triangulations, and there are, in all, $89 - 27 = 62$ non-isomorphic ITPTs obtained by deleting an edge
 from  ${\rm{T}}_i$. As a result, we get Series 3 of ITPTs. This Series is provided in Table II with isomorphic pairs placed
 in one row; as a matter of notation, we write, for instance,  ${\rm{T}}_2-a$ to denote the triangulation obtained from  $\rm{T}_2$ by deleting an arbitrary edge in orbit  $a$.

    \begin{center}
    \hspace*{0.0cm}       \pdfimage width 1.0\textwidth {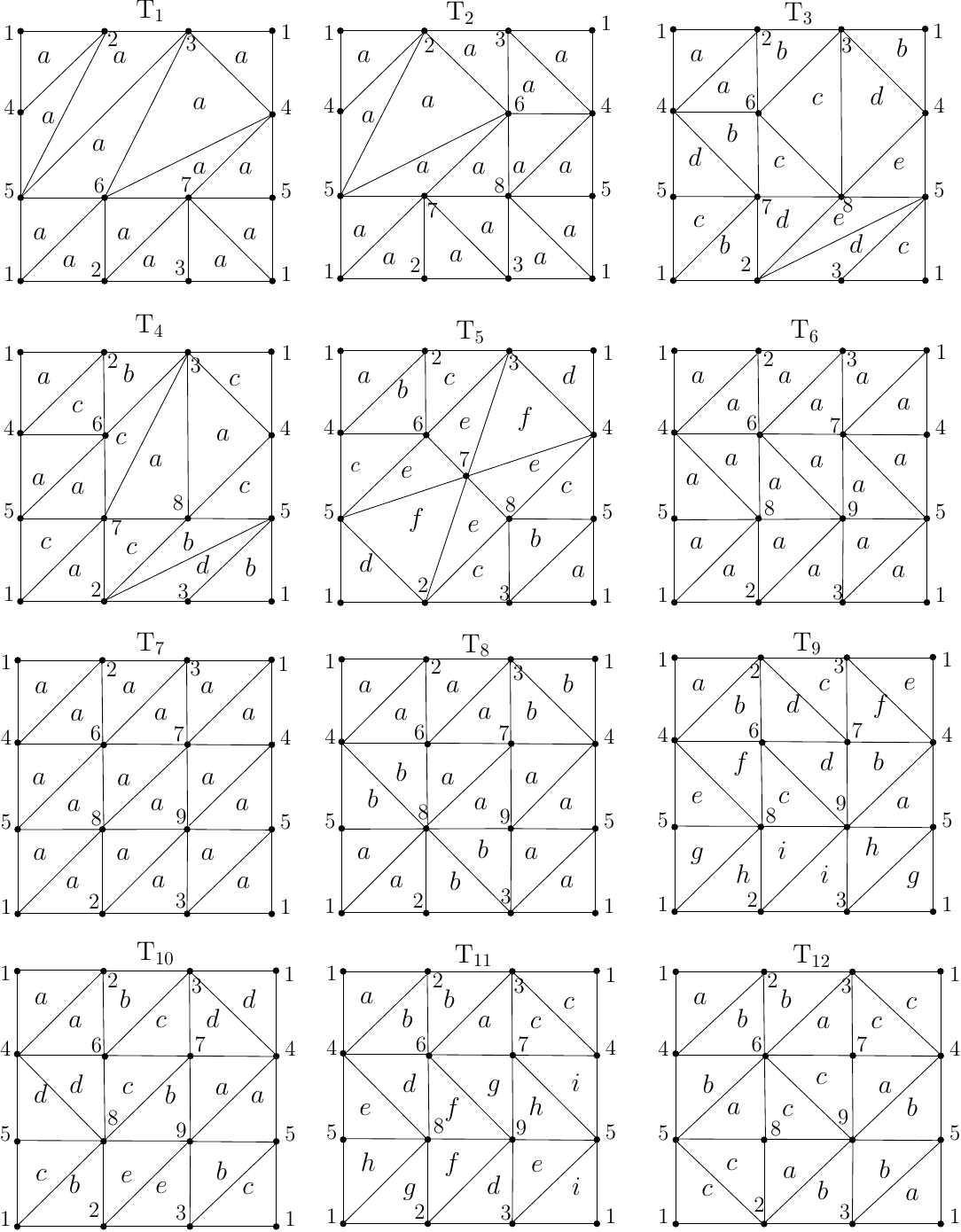}
\newline

    Figure 3. The face orbits (to be continued).
    \end{center}

   \begin{center}
   \hspace*{0.0cm}     \pdfimage width 1.0\textwidth {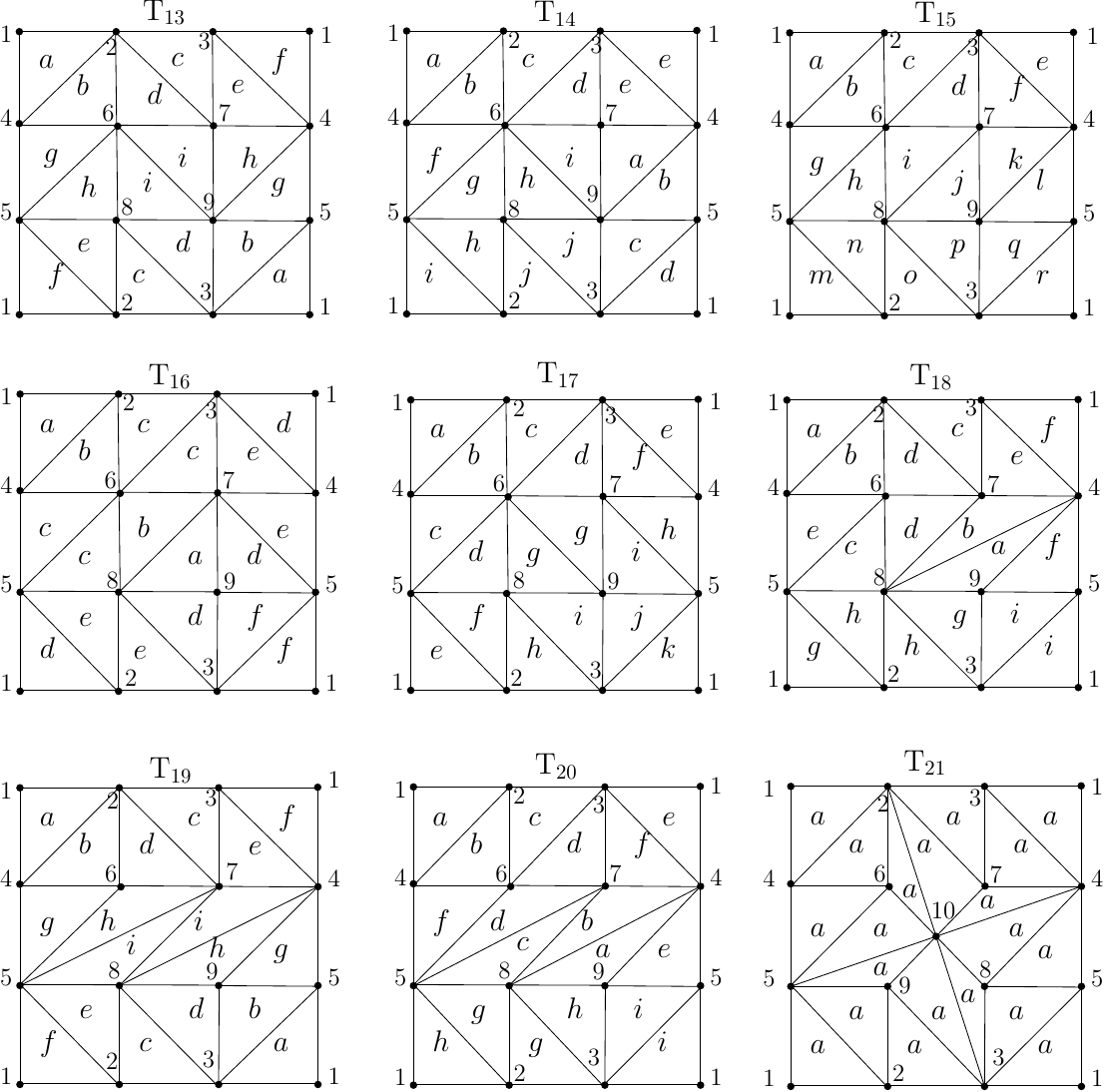}
\newline

    Figure 3. The face orbits (contd.)
    \end{center}

 It remains to check that the $62$ triangulations in Series 3 (Table II) are pairwise non-isomorphic. \ The d-{\it vector} \ ({\it degree vector}) \ of \ a \ triangulation \  $T$  \ is \ defined \ to be  $d(T)=(n_3,n_4, \dots, n_{\mid V(T)\mid -1})$, where   $n_r$ is the number of $r$-valent vertices. The bd-{\it sequence} ({\it boundary degree sequence}) is the cyclic sequence of the degrees of boundary vertices.

 The triangulations in Table II either have differing  d-vectors or   bd-sequences except  for the  following three  non-isomorphic pairs: (i) ${\rm{T}}_2-a\ncong{\rm{T}}_2-b$, (ii) ${\rm{T}}_9-f\ncong{\rm{T}}_{10}-e$, (iii)~${\rm{T}}_9-k\ncong{\rm{T}}_{10}-f$. Proofs of their non-isomorphism are provided in the next couple of paragraphs.

To show that the triangulations in pair (i) are non-isomorphic, we pick  ${\rm{T}}_2-[6,8]$  as a representative of  ${\rm{T}}_2-a$, and pick   ${\rm{T}}_2-[6,7]$  as  ${\rm{T}}_2-b$. Since $6$ and $8$ are the only two $5$-valent vertices in  ${\rm{T}}_2-[6,8]$, and $6$ and $7$ are the only two such vertices in  ${\rm{T}}_2-[6,7]$, and since the edges $[6,8]$ and $[6,7]$ are non-similar in  ${\rm{T}}_2$, it follows that no isomorphism is possible between   ${\rm{T}}_2-[6,8]$  and  ${\rm{T}}_2-[6,7]$.

    \begin{center}
    \hspace*{0.0cm}    \pdfimage width 1.0\textwidth {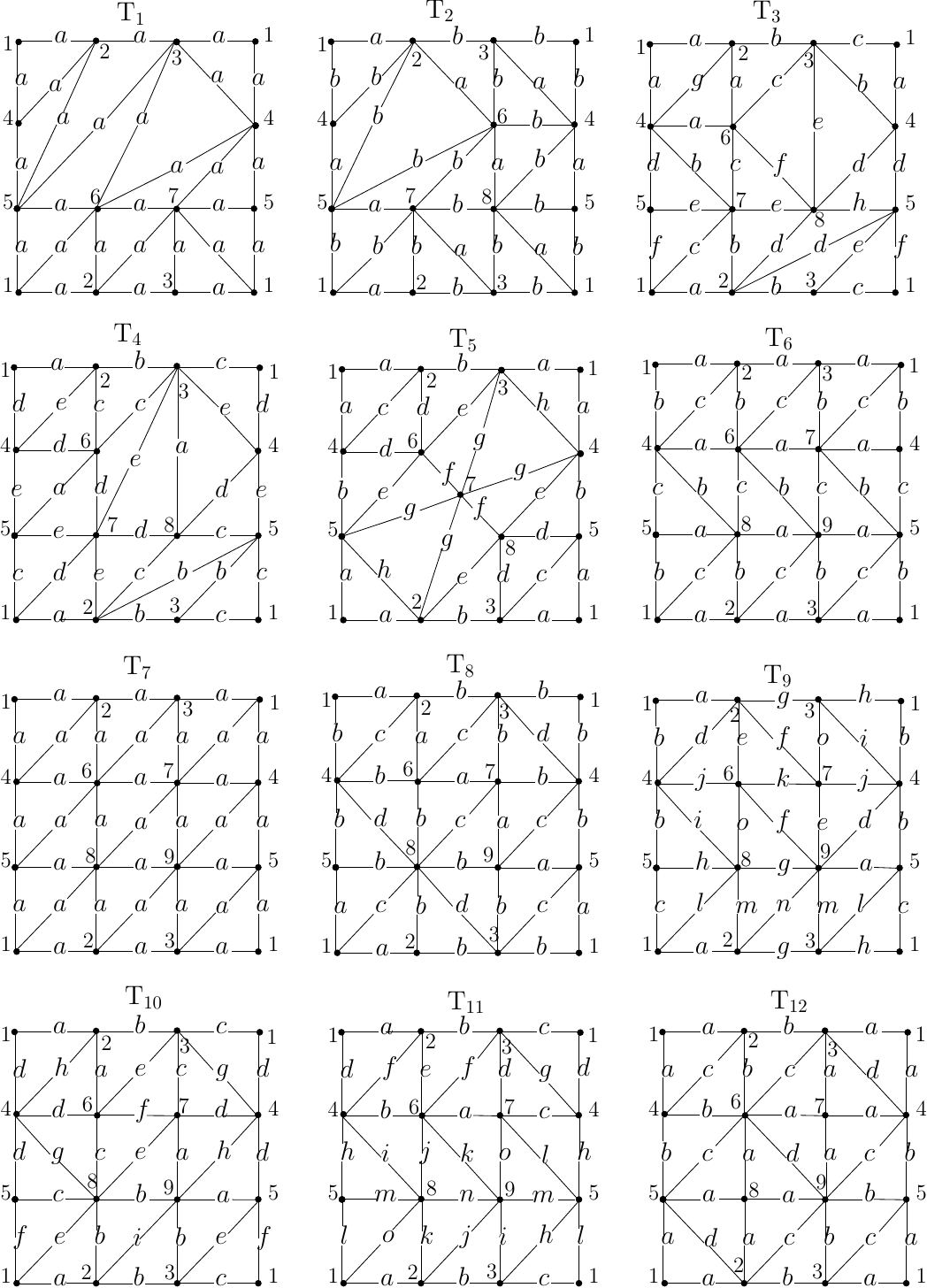}
\newline

    Figure 4. The edge orbits (to be continued).
    \end{center}

The triangulations in pair (ii) are non-isomorphic because only ${\rm{T}}_{10}-e$  has a face with all vertices $6$-valent. Finally, the ones in pair (iii) are non-isomorphic because the only two $6$-valent vertices are adjacent in ${\rm{T}}_{10}-f$  but non-adjacent in  ${\rm{T}}_9-k$.

 \section{The search for ITPTs: Case 3: Identifying non-similar 3-dividers }
\label{S9}

{\it Case 3.}   Parent triangulation  $T^*$  is in $ \Lambda_3$ or can be obtained from a member of   $ \Lambda_3$ by
a sequence of splittings. 
\medskip

    By Eq. (1), if $T^*={\rm{sp}}  \langle  u,v,w \rangle({\rm{T}}_i) \in \Lambda_3$, then  $  \langle u,v,w \rangle $ is a $3$-divider which divides the link of  $v$
    into two edge-disjoint paths, {\it sublinks}, one of which---$u,x,y,w$---has length $3$ and the other has length at least
    $3$; see bottom left of Figure 1. This type of splitting can be thought of as the {\it cracking of the face}
     $[x,v,y]$. In this context, we regard the edge $\varepsilon = [x,y]$  as the {\it base} and the vertex  $v$ as one of the two
     {\it apices} opposite to the base.

\medskip

    \begin{center}
     \hspace*{0cm}      \pdfimage width 1.\textwidth {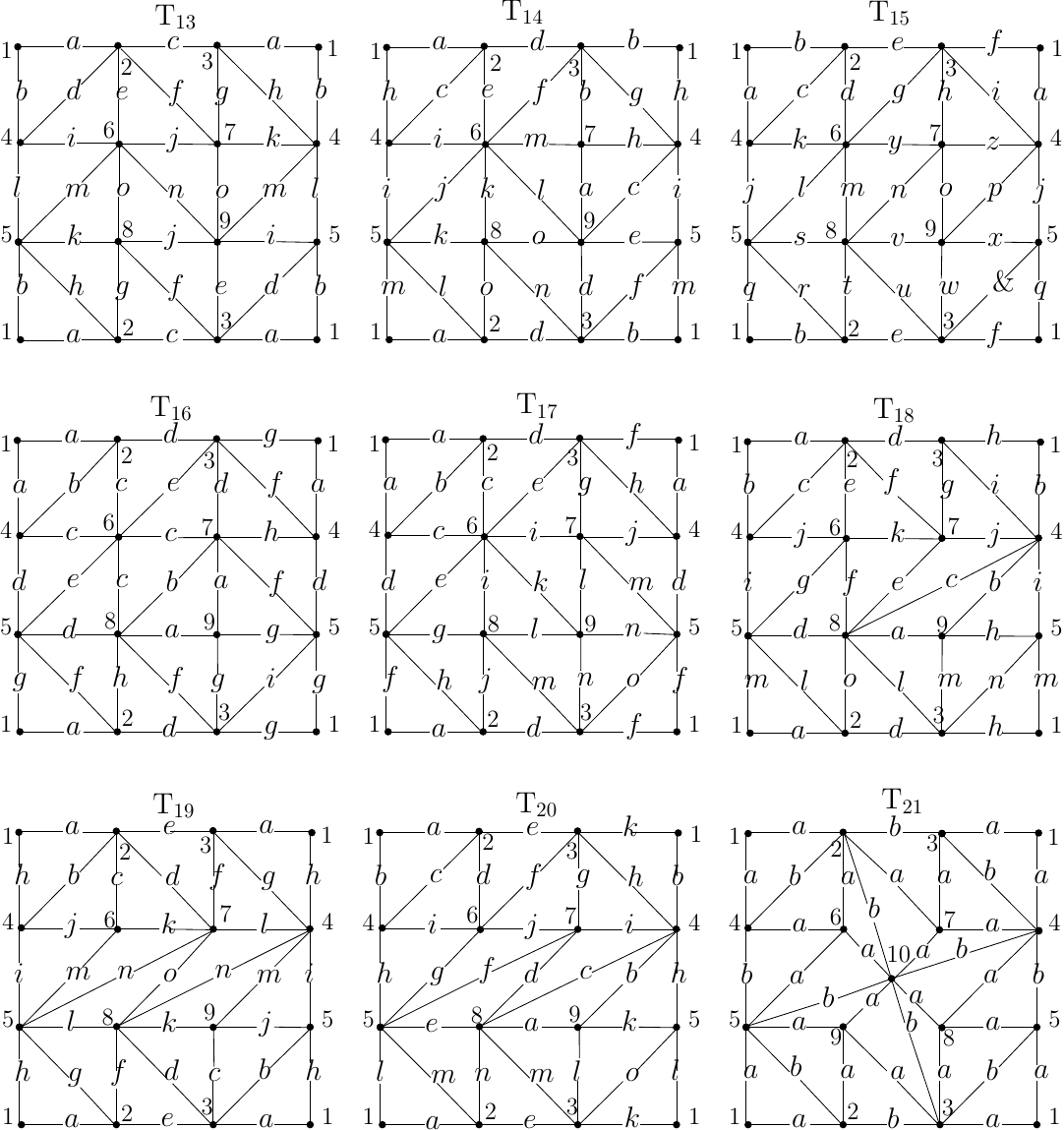}
\newline

    Figure 4. The edge orbits (contd.)
    \end{center}

\newpage
  \begin{center}
 Table II. ITPTs: Series 3 (to be continued).
               \begin{tabular}{|c|l |c|c|}
            \hline  No. & Triangulation      & d-vector & bd-sequence
            \\
            \hline
           1& \rm{T}$_1-a$             &   (0,0,2,5) & (5,6,5,6)\\
            \hline
            2& \rm{T}$_2-a \cong$ \rm{T}$_3-g$  & (0,0,2,6,0) & (5,6,5,6)\\
            \hline
            3& \rm{T}$_2-b$             &  (0,0,2,6,0) & (5,6,5,6) \\
            \hline
           4 & \rm{T}$_3-c$            &  (0,1,2,3,2) & (4,6,5,7)\\
            \hline
            5 & \rm{T}$_3-d$            &  (0,0,3,4,1) & (5,6,6,6)\\
            \hline
           6 & \rm{T}$_3-e$            &  (0,0,4,2,2) & (5,5,5,7)\\
            \hline
            7& \rm{T}$_3-f \cong $\rm{T}$_5-c$ &  (0,1,2,3,2) & (4,6,5,6) \\
            \hline
            8 & \rm{T}$_3-h$            &  (0,0,4,2,2) & (5,7,5,7)\\
            \hline
           9& \rm{T}$_4-a \cong $\rm{T}$_5-h$ &  (0,1,2,3,2) & (4,6,6,6)\\
            \hline
            10 & \rm{T}$_4-b$           &  (0,0,3,4,1) & (5,6,7,6) \\
            \hline
           11 & \rm{T}$_4-d$           &  (0,1,3,1,3) & (4,7,5,7) \\
            \hline
            12 & \rm{T}$_5-d$           &  (0,2,1,2,3) & (4,7,6,7)\\
            \hline
            13 & \rm{T}$_5-e$           &  (0,2,1,2,3) & (4,6,6,7) \\
            \hline
            14 & \rm{T}$_5-f$           &  (0,2,2,0,4) & (4,7,5,7) \\
            \hline
            15 & \rm{T}$_5-g$           &  (0,1,3,1,3) & (5,5,6,7) \\
            \hline
            16& \rm{T}$_6-c \cong $\rm{T}$_{11}-g$ & (0,0,2,7,0,0) & (5,6,5,6) \\
            \hline
            17 & \rm{T}$_8-a$               &  (0,2,4,0,0,3) & (4,8,4,8) \\
            \hline
            18& \rm{T}$_8-c \cong $\rm{T}$_{17}-b$ &  (0,1,5,0,1,2) & (4,5,7,5) \\
            \hline
            19& \rm{T}$_8-d \cong $\rm{T}$_{10}-i$ &  (0,0,6,0,2,1) & (5,7,5,7) \\
            \hline
            20 & \rm{T}$_9-c \cong $\rm{T}$_{18}-n$ &  (0,2,2,2,2,1) & (4,6,4,6) \\
            \hline
            21 & \rm{T}$_9-d \cong $\rm{T}$_{11}-j$ &  (0,0,4,3,2,0) & (5,6,5,7) \\
            \hline
            22 & \rm{T}$_9-f \cong $\rm{T}$_{15}-c$ &  (0,1,3,3,1,1) & (4,5,6,6) \\
            \hline
            23 & \rm{T}$_9-k$                &  (0,2,2,2,2,1) & (4,7,4,7) \\
              \hline
              24 & \rm{T}$_9-l \cong $\rm{T}$_{17}-h$ & (0,1,4,1,2,1) & (4,5,5,7) \\
            \hline
 25 & \rm{T}$_9-n$                &  (0,0,4,4,1) & (6,6,6,6) \\
 \hline

        \end{tabular}
         \end{center}

\medskip

     In this section, we identify all non-similar $3$-dividers in each  ${\rm{T}}_i$ by using inclusion-exclusion. The idea behind this is that instead of counting non-similar $3$-dividers, we judiciously count the base edges   $\varepsilon$  that give rise to them. For this, we associate with each $3$-divider $\langle u,v,w\rangle $  the edge   $\varepsilon = [x,y]$ as indicated on the bottom left of Figure 1 in which the degree of the apex $v$  is assumed to be at least $6$ (in ${\rm{T}}_i$). We say that $\varepsilon$  {\it gives rise to the} $3$-{\it divider} $  \langle u,v,w \rangle $. Each edge (taken as the base) in a triangulation of a closed 2-manifold gives rise to at most two $3$-dividers centered at the apices.

    Let fidx$(\varepsilon)$  denote the  f-{\it index} ({\it face-orbit index}) of an edge  $\varepsilon$  in a given triangulation, defined as follows:  fidx$(\varepsilon)=1$ if the incident faces are in the same orbit (that is, the faces are marked by the same letter in Figure 3),  fidx$(\varepsilon)=2$ if the incident faces are in different orbits (marked by two different letters).

     Let us call an edge $\varepsilon$  of a triangulation $T$  {\it dually reversible} if there is an involutory automorphism of  $T$ that fixes the base  $\varepsilon$ and swaps the two apices (such an automorphism reverses the edge dual to $\varepsilon$). Of course,  $\varepsilon$ is (dually) irreversible whenever  fidx$(\varepsilon)=2$, but may be either reversible or irreversible when  fidx$(\varepsilon)=1$.

     We first study  ${\rm{T}}_i$ for $i=1, 2, \ldots,20$, and postpone  ${\rm{T}}_{21}$ until Lemma \ref{9.5}. By direct inspection (with help from Figures 3, 4 along with Table I), we have verified the following lemma.

\newpage
    \begin{center}
Table II. ITPTs: Series 3 (contd.)
        \begin{tabular}{|c|l |c|c |}
            \hline  No. & Triangulation      & d-vector & bd-sequence
            \\
            \hline
            26 & \rm{T}$_9-o$                &  (0,1,4,1,2,1) & (4,7,5,8) \\
            \hline
            27 & \rm{T}$_{10}-e \cong $\rm{T}$_{15}-r$ &  (0,1,3,3,1,1) & (4,5,6,6) \\
            \hline
            28 & \rm{T}$_{10}-f \cong $\rm{T}$_{20}-o$ &  (0,2,2,2,2,1) & (4,7,4,7)  \\
            \hline
            29 & \rm{T}$_{11}-f$                &  (0,0,3,5,1,0) & (5,5,6,6)  \\
            \hline
            30 & \rm{T}$_{11}-l \cong $\rm{T}$_{15}- \&$ &  (0,1,2,4,2,0) & (4,6,5,7)  \\
            \hline
            31 & \rm{T}$_{11}-o \cong $\rm{T}$_{13}-h$   &  (0,1,2,4,2,0) & (4,6,5,6) \\
            \hline
            32 & \rm{T}$_{12}-c \cong $\rm{T}$_{14}-n$   &  (0,3,0,2,4,0) & (4,6,7,6) \\
            \hline
            33 & \rm{T}$_{12}-d$                  &  (0,3,0,2,4,0) & (4,6,4,6) \\
            \hline
            34 &\rm{T}$_{13}-f \cong $\rm{T}$_{14}-f$ & (0,2,1,3,3,0) & (4,6,6,7) \\
            \hline
            35 & \rm{T}$_{13}-j \cong $\rm{T}$_{18}-c$ & (0,2,2,1,4,0) & (4,6,5,7) \\
            \hline
            36 & \rm{T}$_{13}-m \cong $\rm{T}$_{15}-g$ & (0,1,3,2,3,0) & (5,5,6,7) \\
            \hline
            37 & \rm{T}$_{13}-n$                & (0,1,4,0,4,0) & (5,5,5,5) \\
            \hline
            38 & \rm{T}$_{14}-c \cong $\rm{T}$_{17}-m$ & (0,2,2,2,2,1) & (4,5,7,6)\\
            \hline
            39 & \rm{T}$_{14}-g$                & (0,2,1,3,3,0) & (4,6,4,7)\\
            \hline
            40 & \rm{T}$_{14}-j\cong $\rm{T}$_{16}-e$  & (0,2,1,4,1,1) & (5,6,7,6)\\
            \hline
            41 & \rm{T}$_{14}-l\cong $\rm{T}$_{15}-n$  & (0,2,2,2,2,1) & (4,5,5,6)\\
            \hline
            42 & \rm{T}$_{14}-o$                & (0,3,1,1,3,1) & (4,7,5,8)\\
            \hline
            43 & \rm{T}$_{15}-l\cong $\rm{T}$_{17}-e$  & (0,1,3,3,1,1) & (5,6,6,7)\\
            \hline
            44 & \rm{T}$_{15}-o\cong $\rm{T}$_{19}-b$  & (0,3,0,3,2,1) & (4,6,4,7)\\
            \hline
            45 & \rm{T}$_{15}-p\cong $\rm{T}$_{16}-f$  & (0,2,1,4,1,1) & (4,5,6,7)\\
            \hline
            46 & \rm{T}$_{15}-t$                & (0,1,4,1,2,1) & (5,7,5,8)\\
            \hline
            47 & \rm{T}$_{16}-b\cong $\rm{T}$_{17}-k$  & (0,2,2,3,0,2) & (4,5,6,5)\\
            \hline
            48 & \rm{T}$_{16}-h$                & (0,2,2,3,0,2) & (5,8,5,8)\\
            \hline
            49 & \rm{T}$_{16}-i$                & (0,2,0,5,2,0) & (4,7,4,7)\\
            \hline
            50 & \rm{T}$_{17}-j$                & (0,2,3,1,1,2) & (4,8,5,8)\\
            \hline
            51 & \rm{T}$_{17}-l$                & (0,3,1,2,1,2) & (4,7,4,8)\\
            \hline
            52 & \rm{T}$_{17}-o$                & (0,1,3,2,3,0) & (4,7,5,7)\\
            \hline
            53 & \rm{T}$_{18}-f\cong $\rm{T}$_{19}-n$  & (0,3,1,1,3,1) & (4,5,6,7)\\
            \hline
            54 & \rm{T}$_{18}-g$                & (0,3,1,1,3,1) & (4,7,6,8)\\
            \hline
            55 & \rm{T}$_{18}-k$                & (0,4,0,0,4,1) & (4,7,4,7)\\
            \hline
            56 & \rm{T}$_{18}-o$                & (0,2,2,2,2,1) & (6,7,6,7)\\
            \hline
            57 & \rm{T}$_{19}-d \cong $\rm{T}$_{20}-f$ & (0,3,1,2,1,2) & (4,5,7,6) \\
            \hline
            58 & \rm{T}$_{19}-f$                & (0,3,1,2,1,2) & (5,7,6,8) \\
            \hline
            59 & \rm{T}$_{19}-o$                & (0,3,2,0,2,2) & (5,8,5,8) \\
            \hline
            60 & \rm{T}$_{20}-d$                & (0,3,2,1,0,3) & (4,8,5,8) \\
            \hline
            61 & \rm{T}$_{20}-j$                & (0,4,0,2,0,3) & (4,8,4,8) \\
            \hline
            62 & \rm{T}$_{20}-n$                & (0,2,4,0,0,3) & (5,8,5,8) \\
            \hline
        \end{tabular}
\end{center}

\begin{lemma}\label{9.1}
Each edge of  ${\rm{T}}_i$ $(i=1, 2, \ldots,20)$ with {\rm{f}}-index 1 is dually reversible.
\end{lemma}

\begin{lemma}\label{9.2}
 Let  $i \in \{1, 2, \ldots,20\}$  and let  $\varepsilon \in E({\rm{T}}_i)$. Assume that  $\varepsilon$  gives rise to two 3-dividers. Then the two 3-dividers are similar if and only if   ${\rm{fidx}}(\varepsilon)=1$. (Equivalently, they are non-similar if and only if  ${\rm{fidx}}(\varepsilon)=2$.)
\end{lemma}

\begin{proof}
To prove the ``if'' part, assume  fidx$(\varepsilon)=1$. Then the two $3$-dividers are similar by 
Lemma~\ref{9.1}.

To prove the ``only if'' part, assume that the two $3$-dividers are similar. We have two cases to
consider: (i) the two faces incident with  $\varepsilon$ are similar, and (ii) the two incident faces are
non-similar. In case (i), we immediately come to the desired conclusion:  fidx$(\varepsilon)=1$. In
case~(ii), denote by $v$  and  $z$ the two apices corresponding to  $\varepsilon$. Clearly, any
automorphism that takes the $3$-divider at  $v$ onto the $3$-divider at  $z$  necessarily sends the
sublink of  $v$ that contains  $\varepsilon$ onto the sublink of  $z$  that does not contain
$\varepsilon$, which requires that both sublinks have the same length, so the apices are both
necessarily 6-valent. On the other hand, by direct inspection, we have verified that for each edge
of   ${\rm{T}}_i$ $(i=1, 2, \ldots,20)$ with f-index $2$, the apices are either non-similar or
non-6-valent. Thereby, we still come to the ultimate equality:  fidx$(\varepsilon)=1$.
\end{proof}

     Given a triangulation  $T$, we define the v-{\it index} ({\it vertex-degree index}) of a base
     $\varepsilon_m \in E(T)$  to be the number of its apices with degree at least $6$, and denote it by
     vidx$(\varepsilon_m)$. Clearly,  vidx$(\varepsilon_m)= 0, 1$, or $2$. We define the  $s_1$-{\it invariant } as follows:

\begin{equation} s_1=s_1(T):= \sum_m  {\rm{min}}({\rm{fidx}}(\varepsilon_m), {\rm{vidx}}(\varepsilon_m)),
\end{equation}
where the sum is taken over all edge orbits in  $T$. The general ($m$th) term in Eq. (2) has value
$0$ if and only if  vidx=$0$, has value $2$ if and only if  fidx=vidx=$2$, and has value $1$ in all
other cases. To construct such a function, we use   min(fidx, vidx).

     By Lemma 9.2,  $s_1({\rm{T}}_i)$  is equal to the number of non-similar $3$-dividers in  ${\rm{T}}_i$  $(i =1, 2, \ldots,20)$
     with some of them counted twice, as explained in the next paragraph. The bases $\varepsilon_m$  in Eq. (2) give rise to
     the counted $3$-dividers.

     In the specific case in which the degree of a vertex, $v_n$, is precisely $6$, any two edges opposite to each
     other in the link of  $v_n$, when taken as bases, give rise to the same $3$-divider, centered at  $v_n$, regardless
     of whether the two edges are similar or not. Let kidx$(v_n)$  be the
      $\chi$-{\it index} ({\it link-chirality index}) of a $6$-valent vertex  $v_n \in V(T)$, that is, the
      number of distinct pairs of non-similar opposite edges in the link of  $v_n$ (distinct as unordered pairs of corresponding letters).
      Clearly,  kidx$(v_n)\in \{0,1,2,3 \}$. We define the  $s_2$-{\it invariant} as follows:

 $$ s_2=s_2(T):= \sum_n {\rm{kidx}}(v_n),$$

\noindent where the sum is taken over all $6$-valent vertex orbits in  $T$; if there are no
$6$-valent vertices,
       $s_2:=0$. It is not hard to observe that  $s_2$ is equal to the number of $3$-dividers doubly counted in Eq. (2).
        We finally come to the following simple formula:

\begin{lemma}\label{9.3}
The number of non-similar 3-dividers in  ${\rm{T}}_i \,  ( i =1, 2, \ldots,20)$ is equal to the difference
$s_1({\rm{T}}_i) - s_2({\rm{T}}_i)$.
\end{lemma}

     The following is an example of how to intelligently obtain a complete list of pairwise non-similar
      $3$-dividers in a given triangulation without going through a tedious check. The example addresses the ``hardest triangulation''
      ${\rm{T}}_3$.

\begin{example}\label{9.4}

{\em To illustrate the proof of Lemma \ref{9.3}, which was given prior to its statement, we consider an example of
its use in determining all (pairwise) non-similar $3$-dividers in  ${\rm{T}}_3$. Under the action of  Aut({\rm{T}}$_3$), the vertices of
 ${\rm{T}}_3$ fall into four orbits marked by the letters  $a$--$d$  in Figure 2, the edges fall into eight orbits marked by  $a$--$h$  in Figure 4, and the faces fall into five orbits  $a$--$e$  in Figure 3. We choose, arbitrarily, the following eight edges as representatives of the edge orbits: edge [$2,6]$ in orbit  $a$, $[3,4]$ in  $b$, $[6,7]$ in  $c$, $[4,5]$ in  $d$, $[3,8]$ in  $e$, $[1,5]$ in  $f$, $[2,4]$ in  $g$, $[5,8]$ in  $h$. We calculate the  $s_1$-invariant by writing the edge orbits in alphabetical order and checking with Figure 3---see Table III.

\medskip

    \begin{center}
Table III. Calculation of the $s_1$-invariant.
\newline
\newline
        \begin{tabular}{|l|c|c|c|c|c|c|c|c|c|}
            \hline edge              & [2,6] & [3,4]  & [6,7]  & [4,5]  &[3,8]  & [1,5]  & [2,4] & [5,8] & \\
            \hline orbit             &  $a$   & $b$      &  $ c$    &   $d$    & $e$    &  $f$     &  $ g$  &   $h$   & \\
            \hline f-index         &  2    &   2    &   2    &   2    & 2     &  1     &   1   &   1   & \\
            \hline v-index         &  2    &   1    &   2    &   2    & 1     &  2     &   0   &   2   & \\
            \hline min(fidx,vidx) &  2    &   1    &   2    &   2    & 1     &  1     &   0   &   1   & \makecell{$s_1=\sum$ \\ $=10$}\\
            \hline
        \end{tabular}
    \end{center}

\medskip

$
s_1({\rm{T}}_3)= \Sigma_m \mbox{min} ({\rm{fidx}}(\varepsilon _m), {\rm{vidx}}(\varepsilon _m)) \\ 
 ={\rm{min}}({\rm{fidx}}([2,6]),{\rm{vidx}}( [2,6 ]))+{\rm{min}}({\rm{fidx}}([3,4]),{\rm{vidx}}( [ 3,4 ])) \\ +{\rm{min}}({\rm{fidx}}([6,7]),{\rm{vidx}}( [ 6, 7]))
 +{\rm{min}}({\rm{fidx}}([4,5]),{\rm{vidx}}( [ 4, 5])) \\ +{\rm{min}}({\rm{fidx}}([3,8]),{\rm{vidx}}( [ 3, 8]))+{\rm{min}}({\rm{fidx}}([1,5]),{\rm{vidx}}( [ 1, 5]))\\
 +{\rm{min}}({\rm{fidx}}([2,4]),{\rm{vidx}}( [ 2, 4]))+{\rm{min}}({\rm{fidx}}([5,8]),{\rm{vidx}}( [ 5, 8]) )\\
 ={\rm{min}}(2,2) +{\rm{min}}(2,1)+{\rm{min}}(2,2)+{\rm{min}}(2,2)+{\rm{min}}(2,1) \\+{\rm{min}}(1,2)+{\rm{min}}(1,0)+{\rm{min}}(1,2)\\
 =2+1+2+2+1+1+0+1=10.$

\medskip

Now we calculate the  $s_2$-invariant. There are four orbits into which the vertices of   ${\rm{T}}_3$ fall:  $a,b,c,$ and  $d$; see Figure 2. However, only the vertices in orbits  $c$ and $d$  have degree $6$; we pick vertices $7$ and $8$ as their representatives (respectively). The opposite pairs of edges in the link of vertex $7$ are marked by letters (Figure 4) as follows: $\{a ,a \}$, $\{ f,f \}$, and  $\{d,d \}$; here kidx$(7)=0$  because there is no pair of different letters. The opposite pairs of edges in the link of vertex 8 are:  $\{ c,d \}$, $\{c ,d \}$, and  $\{ b,b \}$; here kidx$(8)=1$  because there is only one pair of different letters, $\{ c, d\}$, distinct from the other pairs. Therefore,

$$ s_2({\rm{T}}_3)= \sum_n {\rm{kidx}}(v_n)={\rm{kidx}}(7)+{\rm{kidx}}(8)=0+1=1. $$

Thus,  ${\rm{T}}_3$ has exactly  $s_1-s_2 =10-1=9$ non-similar $3$-dividers. In fact, we can go even further.
 The proposed approach allows {\it explicitly} identifying all non-similar $3$-dividers in  ${\rm{T}}_3$; see
  below and check with Table III.

Firstly, we inspect the edges with ${\rm{min}}$(fidx,vidx) $=0$ which is equivalent to  vidx $=0$.
 We discard such edges because the $3$-dividers that they give rise to are, in fact, $2$- or $1$-dividers.
 In the example under consideration, we discard the edge $[2,4]$.

Secondly, we inspect the edges with ${\rm{min}}$(fidx, vidx) $=1$ which is equivalent to either condition~(i)
(fidx $=2 )$ {\it AND} (vidx $=1 )$, or (ii) (fidx $=1)$ {\it AND} (vidx $=2)$, or (iii) (fidx $=1)$ {\it AND} (vidx $=1)$.
 In condition (i), since  vidx $=1$, there is only one $3$-divider each such edge gives rise to; in
 the example under consideration, we have two edges in this case: $[3,4]$  and $[3,8]$; the $3$-dividers
 that these edges give rise to are $  \langle  5,8,6 \rangle  $  and  $\langle 1,4,5 \rangle $, respectively (Figure 3).
  In condition (ii), since  fidx $=1$, then there is up to similarity just one $3$-divider that each such edge
   gives rise to, by Lemma 9.2; in the example under consideration, we have two edges in this case: $[1,5]$  and
    $[5,8]$; the $3$-dividers that the edges give rise to are $  \langle 2,7,4 \rangle $  and  $  \langle 7,4,3 \rangle $, respectively (Figure 3).
    In condition (iii), there is only one $3$-divider that each such edge gives rise to, by Lemma 9.2;
    however, in the example under consideration, we have no edges in this case.

Thirdly, among the remaining edges, we inspect the ones with  
$$
{\rm{min}}({\rm{fidx}},{\rm{vidx}})=2
$$
which is
equivalent to condition (fidx $=2$) {\it AND} (vidx $=2)$. By Lemma \ref{9.2}, each such edge gives rise to two
non-similar $3$-dividers---namely, as seen in Figure 3, edge $[2,6]$  gives
rise to $  \langle 1,4,7 \rangle $ and $  \langle 8,3,5 \rangle $     (the latter is similar to  $  \langle 8,7,5 \rangle $---Table I), edge $[6,7]$
gives rise to $  \langle 2,4,5 \rangle $  and $  \langle 3,8,2 \rangle $  (similar to $  \langle 7,8,4 \rangle $), edge [4,5] gives rise to  $  \langle 6,7,1 \rangle $
and $  \langle 3,8,2 \rangle $  (similar to  $  \langle 7,8,4 \rangle $); check with Figure 3. Observe that we
have a duplication: the $3$-divider $  \langle 3,8,2 \rangle $  is doubly counted because on one hand it is given
rise to by two non-similar edges ($[6,7]$  and $[4,5]$), but on the other hand, since vertex $8$ is
$6$-valent, these opposite (in the link of $8$) edges give rise to the \textit{same} $3$-divider---$\langle 3,8,2 \rangle $.}
\end{example}

\newpage
 \begin{center}
 Table IV. Calculation of the number of 3-dividers.
\newline
\newline
    \scalebox{1.0}{
        \begin{tabular}{|c|c|c|c|}
            \hline T$_i$  & $s_1$ & $s_2$  & The $\#$ of \\
            &       &        & 3-dividers  \\
            \hline T$_1$      & 1     &  0     &  1          \\
            \hline T$_2$      & 2     &  0     &  2          \\
            \hline T$_3$      & 10    &  1     &  9          \\
            \hline T$_4$      & 6     &  1     &  5          \\
            \hline T$_5$      & 9     &  0     &  9          \\
            \hline T$_6$      & 3     &  0     &  3          \\
            \hline T$_7$      & 1     &  0     &  1          \\
            \hline T$_8$      & 2     &  0     &  2          \\
            \hline T$_9$      & 17    &  3     &  14         \\
            \hline T$_{10}$     & 9     &  1     &  8          \\
            \hline T$_{11}$     & 22    &  6     &  16         \\
            \hline T$_{12}$     & 4     &  0     &  4          \\
            \hline T$_{13}$     & 20    &  3     &  17         \\
            \hline T$_{14}$     & 21    &  3     &  18         \\
            \hline T$_{15}$     & 40    &  9     &  31         \\
            \hline T$_{16}$     & 12    &  3     &  9          \\
            \hline T$_{17}$     & 18    &  3     &  15         \\
            \hline T$_{18}$     & 18    &  0     &  18         \\
            \hline T$_{19}$     & 21    &  3     &  18         \\
            \hline T$_{20}$     & 18    &  3     &  15         \\
            \hline T$_{21}$     & 1     &  0      &  2$^*$          \\
            \hline        &       &        &  $\sum=217$ \\
            \hline
        \end{tabular}
}
       \end{center}
\medskip

     It is easy to count $s_1-s_2$  with help from Figures 3, 4.
     The results are collected in Table IV; this table can be regarded as a corollary of Lemma \ref{9.3}; the asterisk in the last row indicates that the formula in Lemma \ref{9.3} does not apply for ${\rm{T}}_{21}$ addressed in Lemma \ref{9.5}. Moreover,
      similarly to Example \ref{9.4}, we have produced a complete list of non-similar $3$-dividers in each
       ${\rm{T}}_i $ for $  i =1, 2, \ldots, 20$; there are totally $215$ $3$-dividers. Check with Table IV. In addition,
        two more are found in  ${\rm{T}}_{21}$:

\begin{lemma}
\label{9.5} If an edge   $\varepsilon$ of  ${\rm{T}}_{21}$ gives rise to at least one 3-divider, then
$\varepsilon$  is similar to the edge $[2,6]$ (Figure 4) which, in fact, gives
rise to two non-similar 3-dividers: $\langle 5,10,7 \rangle $   and $\langle 5,4,1 \rangle $ (the latter is similar to $\langle 5,10,8 \rangle $---Table I).
\end{lemma}

\begin{proof}
All edges of  ${\rm{T}}_{21}$  split into two orbits---$a$  and $b$  (see Figure 4). For each edge in orbit
 $b$  as the base, both apices are $4$-valent. Therefore, no edge in  $b$  gives rise to a 3-divider. On the other hand, each
 edge in orbit   $a$ is dually irreversible (check with Figure 4 and Table I) and thereby gives rise to two
 non-similar $3$-dividers even though it has f-index $1$ (Figure 3).
\end{proof}

\begin{center}

Table V. $\Lambda_3$ with isomorphic duplications (to be continued). 
		\scalebox{0.85}{
			\begin{tabular}{|ccccc|c|}
				\hline                          &                          &           Triangulations              &                         &                    & Total  \\
				&                         &                         &                         &                    &  \\
				\hline  \fbox{${\rm{sp}} \langle1,7,6 \rangle (  \rm{T}_1)$} &                         &                         &                         &                    &   1     \\
				\hline  ${\rm{sp}} \langle7,8,5 \rangle (  \rm{T}_2)$        & \fbox{${\rm{sp}} \langle6,8,1 \rangle (  \rm{T}_2)$} &                         &                         &                    &   2     \\
				\hline ${\rm{sp}} \langle1,4,7 \rangle (  \rm{T}_3)$         &  ${\rm{sp}} \langle1,4,5 \rangle (  \rm{T}_3)$       &  ${\rm{sp}} \langle2,4,5 \rangle (  \rm{T}_3)$       & \fbox{${\rm{sp}} \langle7,4,3 \rangle (  \rm{T}_3)$} &  ${\rm{sp}} \langle6,7,1 \rangle (  \rm{T}_3)$  &   \\
				${\rm{sp}} \langle8,7,5 \rangle (  \rm{T}_3)$         & \fbox{${\rm{sp}} \langle2,7,4 \rangle (  \rm{T}_3)$} & \fbox{${\rm{sp}} \langle7,8,4 \rangle (  \rm{T}_3)$} & \fbox{${\rm{sp}} \langle5,8,6 \rangle (  \rm{T}_3)$} &                    &  9 \\
				\hline ${\rm{sp}} \langle4,3,2 \rangle (  \rm{T}_4)$         &  ${\rm{sp}} \langle4,3,6 \rangle (  \rm{T}_4)$       &  ${\rm{sp}} \langle8,3,2 \rangle (  \rm{T}_4)$       & \fbox{${\rm{sp}} \langle6,3,1 \rangle (  \rm{T}_4)$} &   ${\rm{sp}} \langle8,7,5 \rangle (  \rm{T}_4)$ &  5 \\
				\hline ${\rm{sp}} \langle6,2,5 \rangle (  \rm{T}_5)$         &  ${\rm{sp}} \langle6,2,7 \rangle (  \rm{T}_5)$       &\fbox{${\rm{sp}} \langle3,2,1 \rangle (  \rm{T}_5)$}  & \fbox{${\rm{sp}} \langle3,2,5 \rangle (  \rm{T}_5)$} & ${\rm{sp}} \langle8,2,4 \rangle (  \rm{T}_5)$   &    \\
				${\rm{sp}} \langle8,2,1 \rangle (  \rm{T}_5)$          &  ${\rm{sp}} \langle7,2,4 \rangle (  \rm{T}_5)$       &  \fbox{${\rm{sp}} \langle3,7,2 \rangle (  \rm{T}_5)$}& \fbox{${\rm{sp}} \langle8,7,6 \rangle (  \rm{T}_5)$} &                    &  9 \\
				\hline ${\rm{sp}} \langle8,6,3 \rangle (  \rm{T}_6)$         & ${\rm{sp}} \langle7,6,4 \rangle (  \rm{T}_6)$        & ${\rm{sp}} \langle9,6,2 \rangle (  \rm{T}_6)$        &                         &                    &  3 \\
				\hline ${\rm{sp}} \langle4,6,7 \rangle (  \rm{T}_7)$         &                         &                         &                         &                    &  1 \\
				\hline ${\rm{sp}} \langle1,8,9 \rangle (  \rm{T}_8)$         & ${\rm{sp}} \langle2,8,4 \rangle (  \rm{T}_8)$        &                         &                         &                    &  2 \\
				\hline ${\rm{sp}} \langle2,9,4 \rangle (  \rm{T}_9)$         &  ${\rm{sp}} \langle2,9,7 \rangle (  \rm{T}_9)$       &  ${\rm{sp}} \langle3,9,6 \rangle (  \rm{T}_9)$       & ${\rm{sp}} \langle3,9,7 \rangle (  \rm{T}_9)$        & ${\rm{sp}} \langle5,9,8 \rangle (  \rm{T}_9)$   &    \\
				${\rm{sp}} \langle5,9,6 \rangle (  \rm{T}_9)$        &  ${\rm{sp}} \langle4,9,8 \rangle (  \rm{T}_9)$       &  ${\rm{sp}} \langle1,8,6 \rangle (  \rm{T}_9)$       &  ${\rm{sp}} \langle5,8,9 \rangle (  \rm{T}_9)$       & ${\rm{sp}} \langle4,8,2 \rangle (  \rm{T}_9)$   &    \\
				${\rm{sp}} \langle3,4,5 \rangle (  \rm{T}_9)$        &  ${\rm{sp}} \langle5,4,2 \rangle (  \rm{T}_9)$       & ${\rm{sp}} \langle2,4,7 \rangle (  \rm{T}_9)$        & ${\rm{sp}} \langle6,4,3 \rangle (  \rm{T}_9)$        &                    & 14 \\
				\hline  ${\rm{sp}} \langle5,4,2 \rangle (  \rm{T}_{10})$     &  ${\rm{sp}} \langle5,4,3 \rangle (  \rm{T}_{10})$    &  ${\rm{sp}} \langle8,9,5 \rangle (  \rm{T}_{10})$    & ${\rm{sp}} \langle2,9,4 \rangle (  \rm{T}_{10})$     &${\rm{sp}} \langle7,8,5 \rangle (  \rm{T}_{10})$ &    \\
				${\rm{sp}} \langle1,8,7 \rangle (  \rm{T}_{10})$     &  ${\rm{sp}} \langle9,8,4 \rangle (  \rm{T}_{10})$    &  ${\rm{sp}} \langle9,8,5 \rangle (  \rm{T}_{10})$    &                         &                    &  8  \\
				\hline  ${\rm{sp}} \langle3,6,8 \rangle (  \rm{T}_{11})$     & ${\rm{sp}} \langle7,6,4 \rangle (  \rm{T}_{11})$     &  ${\rm{sp}} \langle9,6,2 \rangle (  \rm{T}_{11})$    & ${\rm{sp}} \langle7,3,5 \rangle (  \rm{T}_{11})$     & ${\rm{sp}} \langle7,3,9 \rangle (  \rm{T}_{11})$&    \\
				${\rm{sp}} \langle6,3,5 \rangle (  \rm{T}_{11})$      &  ${\rm{sp}} \langle6,3,1 \rangle (  \rm{T}_{11})$    &  ${\rm{sp}} \langle2,3,1 \rangle (  \rm{T}_{11})$    & ${\rm{sp}} \langle2,3,4 \rangle (  \rm{T}_{11})$     & ${\rm{sp}} \langle9,3,4 \rangle (  \rm{T}_{11})$&    \\
				${\rm{sp}} \langle4,5,3 \rangle (  \rm{T}_{11})$      & ${\rm{sp}} \langle8,5,9 \rangle (  \rm{T}_{11})$     & ${\rm{sp}} \langle1,5,7 \rangle (  \rm{T}_{11})$     & ${\rm{sp}} \langle6,8,1 \rangle (  \rm{T}_{11})$     & ${\rm{sp}} \langle9,8,5 \rangle (  \rm{T}_{11})$&    \\
				${\rm{sp}} \langle2,8,4 \rangle (  \rm{T}_{11})$     &                         &                         &                         &                    & 16 \\
				\hline ${\rm{sp}} \langle9,6,4 \rangle (  \rm{T}_{12})$      & ${\rm{sp}} \langle4,6,7 \rangle (  \rm{T}_{12})$     & ${\rm{sp}} \langle8,6,3 \rangle (  \rm{T}_{12})$     &  ${\rm{sp}} \langle5,6,3 \rangle (  \rm{T}_{12})$    &                    & 4  \\
				\hline  ${\rm{sp}} \langle1,2,7 \rangle (  \rm{T}_{13})$     &  ${\rm{sp}} \langle1,2,3 \rangle (  \rm{T}_{13})$    & ${\rm{sp}} \langle5,2,7 \rangle (  \rm{T}_{13})$     & ${\rm{sp}} \langle5,2,6 \rangle (  \rm{T}_{13})$     & ${\rm{sp}} \langle8,2,4 \rangle (  \rm{T}_{13})$&    \\
				${\rm{sp}} \langle8,2,6 \rangle (  \rm{T}_{13})$     & ${\rm{sp}} \langle 3,2,4 \rangle (  \rm{T}_{13})$     & ${\rm{sp}} \langle 9,4,1 \rangle (  \rm{T}_{13})$     & ${\rm{sp}} \langle 9,4,2 \rangle (  \rm{T}_{13})$     &${\rm{sp}} \langle 7,4,2 \rangle (  \rm{T}_{13})$ &    \\
				${\rm{sp}} \langle 7,4,6 \rangle (  \rm{T}_{13})$     &  ${\rm{sp}} \langle 3,4,6 \rangle (  \rm{T}_{13})$    &  ${\rm{sp}} \langle 3,4,5 \rangle (\rm{T}_{13})$    &  ${\rm{sp}} \langle 1,4,5 \rangle (\rm{T}_{13})$    &${\rm{sp}} \langle 5,6,7 \rangle (\rm{T}_{13})$ &    \\
				${\rm{sp}} \langle 4,6,9 \rangle (\rm{T}_{13})$    &  ${\rm{sp}} \langle 2,6,8 \rangle (\rm{T}_{13})$    &                         &                         &                    & 17 \\
				\hline  ${\rm{sp}} \langle 5,6,3 \rangle (\rm{T}_{14})$     &  ${\rm{sp}} \langle 5,6,7 \rangle (\rm{T}_{14})$    & ${\rm{sp}} \langle 4,6,7 \rangle (\rm{T}_{14})$     &  ${\rm{sp}} \langle 4,6,9 \rangle (\rm{T}_{14})$    &${\rm{sp}} \langle 2,6,8 \rangle (\rm{T}_{14})$ &    \\
				${\rm{sp}} \langle 2,6,9 \rangle (\rm{T}_{14})$     & ${\rm{sp}} \langle 3,6,8 \rangle (\rm{T}_{14})$     & ${\rm{sp}} \langle 2,3,5 \rangle (\rm{T}_{14})$     & ${\rm{sp}} \langle 2,3,4 \rangle (\rm{T}_{14})$     &${\rm{sp}} \langle 8,3,1 \rangle (\rm{T}_{14})$ &    \\
				${\rm{sp}} \langle 1,3,6 \rangle (\rm{T}_{14})$     & ${\rm{sp}} \langle 3,9,7 \rangle (\rm{T}_{14})$     & ${\rm{sp}} \langle 8,9,4 \rangle (\rm{T}_{14})$     & ${\rm{sp}} \langle 6,9,5 \rangle (\rm{T}_{14})$     &${\rm{sp}} \langle 9,4,1 \rangle (\rm{T}_{14})$ &    \\
				${\rm{sp}} \langle 9,4,2 \rangle (\rm{T}_{14})$     & ${\rm{sp}} \langle 3,4,5 \rangle (\rm{T}_{14})$     & ${\rm{sp}} \langle 1,4,5 \rangle (\rm{T}_{14})$     &                         &                    & 18 \\
				\hline  ${\rm{sp}} \langle 1,2,3 \rangle (\rm{T}_{15})$     &  ${\rm{sp}} \langle 5,2,6 \rangle (\rm{T}_{15})$    &  ${\rm{sp}} \langle 8,2,4 \rangle (\rm{T}_{15})$    &  ${\rm{sp}} \langle 9,4,1 \rangle (\rm{T}_{15})$    & ${\rm{sp}} \langle 9,4,2 \rangle (\rm{T}_{15})$&    \\
				${\rm{sp}} \langle 7,4,2 \rangle (\rm{T}_{15})$     &  ${\rm{sp}} \langle 7,4,6 \rangle (\rm{T}_{15})$    &  ${\rm{sp}} \langle 3,4,6 \rangle (\rm{T}_{15})$    &  ${\rm{sp}} \langle 3,4,5 \rangle (\rm{T}_{15})$    &${\rm{sp}} \langle 1,4,5 \rangle (\rm{T}_{15})$ &    \\
				${\rm{sp}} \langle 6,8,3 \rangle (\rm{T}_{15})$    &  ${\rm{sp}} \langle 7,8,2 \rangle (\rm{T}_{15})$    &  ${\rm{sp}} \langle 9,8,5 \rangle (\rm{T}_{15})$    &  ${\rm{sp}} \langle 9,3,4 \rangle (\rm{T}_{15})$    & ${\rm{sp}} \langle 9,3,6 \rangle (\rm{T}_{15})$&    \\
				${\rm{sp}} \langle 5,3,7 \rangle (\rm{T}_{15})$     & ${\rm{sp}} \langle 5,3,2 \rangle (\rm{T}_{15})$     & ${\rm{sp}} \langle 1,3,8 \rangle (\rm{T}_{15})$     & ${\rm{sp}} \langle 1,3,6 \rangle (\rm{T}_{15})$     & ${\rm{sp}} \langle 4,3,2 \rangle (\rm{T}_{15})$&    \\
				${\rm{sp}} \langle 7,3,8 \rangle (\rm{T}_{15})$     & ${\rm{sp}} \langle 2,6,8 \rangle (\rm{T}_{15})$     & ${\rm{sp}} \langle 3,6,5 \rangle (\rm{T}_{15})$     & ${\rm{sp}} \langle 7,6,4 \rangle (\rm{T}_{15})$     & ${\rm{sp}} \langle 9,5,8 \rangle (\rm{T}_{15})$&    \\
				${\rm{sp}} \langle 9,5,2 \rangle (\rm{T}_{15})$     & ${\rm{sp}} \langle 4,5,1 \rangle (\rm{T}_{15})$     & ${\rm{sp}} \langle 4,5,2 \rangle (\rm{T}_{15})$     & ${\rm{sp}} \langle 6,5,3 \rangle (\rm{T}_{15})$     &${\rm{sp}} \langle 6,5,1 \rangle (\rm{T}_{15})$ &     \\
				${\rm{sp}} \langle 8,5,3 \rangle (\rm{T}_{15})$     &                         &                         &                         &                    & 31 \\
				\hline ${\rm{sp}} \langle 5,8,9 \rangle (\rm{T}_{16})$       & ${\rm{sp}} \langle 6,8,3 \rangle (\rm{T}_{16})$      &  ${\rm{sp}} \langle 7,8,2 \rangle (\rm{T}_{16})$     &  ${\rm{sp}} \langle 4,6,7 \rangle (\rm{T}_{16})$     &${\rm{sp}} \langle 3,6,5 \rangle (\rm{T}_{16})$  &    \\
				${\rm{sp}} \langle 2,3,5 \rangle (\rm{T}_{16})$       & ${\rm{sp}} \langle 2,3,4 \rangle (\rm{T}_{16})$      & ${\rm{sp}} \langle 8,3,1 \rangle (\rm{T}_{16})$      & ${\rm{sp}} \langle 1,3,6 \rangle (\rm{T}_{16})$      &                    & 9  \\
				\hline  ${\rm{sp}} \langle 9,3,4 \rangle (\rm{T}_{17})$      & ${\rm{sp}} \langle 9,3,6 \rangle (\rm{T}_{17})$      & ${\rm{sp}} \langle 5,3,2 \rangle (\rm{T}_{17})$      &  ${\rm{sp}} \langle 5,3,7 \rangle (\rm{T}_{17})$     & ${\rm{sp}} \langle 1,3,8 \rangle (\rm{T}_{17})$ &    \\
				${\rm{sp}} \langle 1,3,6 \rangle (\rm{T}_{17})$      & ${\rm{sp}} \langle 4,3,2 \rangle (\rm{T}_{17})$      &  ${\rm{sp}} \langle 7,3,8 \rangle (\rm{T}_{17})$     & ${\rm{sp}} \langle 1,2,3 \rangle (\rm{T}_{17})$      & ${\rm{sp}} \langle 5,2,6 \rangle (\rm{T}_{17})$ &    \\
				${\rm{sp}} \langle 8,2,4 \rangle (\rm{T}_{17})$      & ${\rm{sp}} \langle 2,6,9 \rangle (\rm{T}_{17})$      &  ${\rm{sp}} \langle 3,6,5 \rangle (\rm{T}_{17})$     &  ${\rm{sp}} \langle 7,6,4 \rangle (\rm{T}_{17})$     & ${\rm{sp}} \langle 7,6,5 \rangle (\rm{T}_{17})$ & 15 \\
				\hline
			\end{tabular}
		}   
\end{center}
\medskip

\newpage

		\begin{center}
Table V. $\Lambda_3$ with isomorphic duplications (contnd.)
\newline
\newline
			\scalebox{0.85}{
				\begin{tabular}{|ccccc|c|}
					\hline                           &                         &       Triangulations              &                         &                    &  Total  \\
					&                         &                         &                         &                    &  \\
					\hline  ${\rm{sp}} \langle 2,8,7 \rangle (\rm{T}_{18})$      & ${\rm{sp}} \langle 2,8,4 \rangle (\rm{T}_{18})$      &  ${\rm{sp}} \langle 5,8,4 \rangle (\rm{T}_{18})$     &  ${\rm{sp}} \langle 5,8,9 \rangle (\rm{T}_{18})$     & ${\rm{sp}} \langle 6,8,3 \rangle (\rm{T}_{18})$ &    \\
					${\rm{sp}} \langle 6,8,9 \rangle (\rm{T}_{18})$      & ${\rm{sp}} \langle 7,8,3 \rangle (\rm{T}_{18})$      &  ${\rm{sp}} \langle 2,3,5 \rangle (\rm{T}_{18})$     &   ${\rm{sp}} \langle 2,3,1 \rangle (\rm{T}_{18})$    & ${\rm{sp}} \langle 8,3,1 \rangle (\rm{T}_{18})$ &    \\
					${\rm{sp}} \langle 8,3,4 \rangle (\rm{T}_{18})$      & ${\rm{sp}} \langle 9,3,4 \rangle (\rm{T}_{18})$      & ${\rm{sp}} \langle 9,3,7 \rangle (\rm{T}_{18})$      & ${\rm{sp}} \langle 5,3,7 \rangle (\rm{T}_{18})$      & ${\rm{sp}} \langle 7,4,2 \rangle (\rm{T}_{18})$ &    \\
					${\rm{sp}} \langle 3,4,6 \rangle (\rm{T}_{18})$      & ${\rm{sp}} \langle 3,4,9 \rangle (\rm{T}_{18})$      &  ${\rm{sp}} \langle 2,4,9 \rangle (\rm{T}_{18})$     &                          &                     & 18 \\
					\hline  ${\rm{sp}} \langle 5,7,3 \rangle (\rm{T}_{19})$      & ${\rm{sp}} \langle 6,7,4 \rangle (\rm{T}_{19})$      &  ${\rm{sp}} \langle 2,7,8 \rangle (\rm{T}_{19})$     & ${\rm{sp}} \langle 4,2,8 \rangle (\rm{T}_{19})$      &${\rm{sp}} \langle 4,2,3 \rangle (\rm{T}_{19})$  &    \\
					${\rm{sp}} \langle 1,2,3 \rangle (\rm{T}_{19})$      & ${\rm{sp}} \langle 1,2,7 \rangle (\rm{T}_{19})$      &  ${\rm{sp}} \langle 5,2,6 \rangle (\rm{T}_{19})$     & ${\rm{sp}} \langle 5,2,7 \rangle (\rm{T}_{19})$      &${\rm{sp}} \langle 8,2,6 \rangle (\rm{T}_{19})$  &    \\
					${\rm{sp}} \langle 1,4,8 \rangle (\rm{T}_{19})$      & ${\rm{sp}} \langle 1,4,5 \rangle (\rm{T}_{19})$      &  ${\rm{sp}} \langle 2,4,7 \rangle (\rm{T}_{19})$     & ${\rm{sp}} \langle 2,4,9 \rangle (\rm{T}_{19})$      &${\rm{sp}} \langle 6,4,3 \rangle (\rm{T}_{19})$  &    \\
					${\rm{sp}} \langle 6,4,8 \rangle (\rm{T}_{19})$      & ${\rm{sp}} \langle 5,4,7 \rangle (\rm{T}_{19})$      & ${\rm{sp}} \langle 9,4,3 \rangle (\rm{T}_{19})$      &                          &                     & 18 \\
					\hline ${\rm{sp}} \langle 5,8,9 \rangle (\rm{T}_{20})$       & ${\rm{sp}} \langle 7,8,3 \rangle (\rm{T}_{20})$      & ${\rm{sp}} \langle 4,8,2 \rangle (\rm{T}_{20})$      & ${\rm{sp}} \langle 2,3,5 \rangle (\rm{T}_{20})$      &${\rm{sp}} \langle 2,3,4 \rangle (\rm{T}_{20})$  &     \\
					${\rm{sp}} \langle 8,3,1 \rangle (\rm{T}_{20})$       & ${\rm{sp}} \langle 8,3,7 \rangle (\rm{T}_{20})$      & ${\rm{sp}} \langle 9,3,6 \rangle (\rm{T}_{20})$      &  ${\rm{sp}} \langle 9,3,4 \rangle (\rm{T}_{20})$     &${\rm{sp}} \langle 5,3,7 \rangle (\rm{T}_{20})$  &     \\
					${\rm{sp}} \langle 1,3,6 \rangle (\rm{T}_{20})$       &${\rm{sp}} \langle 7,4,5 \rangle (\rm{T}_{20})$       & ${\rm{sp}} \langle 3,4,9 \rangle (\rm{T}_{20})$      & ${\rm{sp}} \langle 1,4,8 \rangle (\rm{T}_{20})$      & ${\rm{sp}} \langle 6,4,8 \rangle (\rm{T}_{20})$ & 15\\
					\hline  ${\rm{sp}} \langle 5,10,7 \rangle (\rm{T}_{21})$     & ${\rm{sp}} \langle 5,10,8 \rangle (\rm{T}_{21})$     &                          &                          &                     & 2\\
					\hline
				\end{tabular}
			}
		\end{center}
\medskip

\medskip

     Using Lemma \ref{9.3} (as in Example \ref{9.4}) and Lemma \ref{9.5}, we collect in Table V the $217$ triangulations
     produced by splitting all the $217$ non-similar $3$-dividers in ${\rm{T}}_i$ $ (i =1, 2, \ldots,21)$, which provides the whole set $\Lambda_3$  (with isomorphic duplications).

\section{The search for ITPTs: Case 3: Series 4-6}
\label{S10}
In order to apply Corollary \ref{4.2}, we identify $\ast$-triangulations from the $217$ triangulations in Table V.
 Moreover, we retain only non-isomorphic triangulations. It is a matter of mere inspection to obtain the following lemma.

\begin{lemma}
\label{10.1}
There are precisely eleven $\ast$-triangulations in  $\Lambda_3$, as follows: \\
${\rm{sp}} \langle 1,7,6 \rangle ({\rm{T}}_1), {\rm{sp}}\langle  6,8,1 \rangle ({\rm{T}}_2), {\rm{sp}}\langle  7,4,3 \rangle ({\rm{T}}_3),
{\rm{sp}}\langle  2,7,4 \rangle ({\rm{T}}_3), {\rm{sp}}\langle  7,8,4 \rangle ({\rm{T}}_3)$, \\
${\rm{sp}} \langle 5,8,6 \rangle ({\rm{T}}_3), {\rm{sp}}\langle  6,3,1 \rangle ({\rm{T}}_4), {\rm{sp}}\langle  3,2,1 \rangle ({\rm{T}}_5), {\rm{sp}}\langle  3,2,5 \rangle ({\rm{T}}_5), {\rm{sp}}\langle  3,7,2 \rangle ({\rm{T}}_5)$,  \\ and
${\rm{sp}} \langle 8,7,6 \rangle ({\rm{T}}_5)$.

These triangulations are boxed in Table V and are also presented in Figure 5 with all ropes in bold type.
\end{lemma}

It should be noticed that eight of the eleven $\ast$-triangulations in Lemma \ref{10.1} have only one rope
and two $\ast$-vertices. The removal of a $\ast$-vertex in each of the eleven triangulations creates the following $19$ ITPTs: \\

\noindent ${\rm{sp}} \langle 1,7,6 \rangle ({\rm{T}}_1)-7'$, ${\rm{sp}} \langle 1,7,6 \rangle ({\rm{T}}_1)-7''$, ${\rm{sp}} \langle 6,8,1 \rangle ({\rm{T}}_2)-8'$, \\
${\rm{sp}} \langle 6,8,1 \rangle ({\rm{T}}_2)-8''$,  ${\rm{sp}} \langle 7,4,3 \rangle ({\rm{T}}_3)-4'$,  ${\rm{sp}} \langle 7,4,3 \rangle ({\rm{T}}_3)-4''$, \\
${\rm{sp}} \langle 2,7,4 \rangle ({\rm{T}}_3)-7'$, ${\rm{sp}} \langle 2,7,4 \rangle ({\rm{T}}_3)-7''$, ${\rm{sp}} \langle 7,8,4 \rangle ({\rm{T}}_3)-8'$, \\
${\rm{sp}} \langle 5,8,6 \rangle ({\rm{T}}_3)-8'$, ${\rm{sp}} \langle 5,8,6 \rangle ({\rm{T}}_3)-8''$,  ${\rm{sp}} \langle 6,3,1 \rangle ({\rm{T}}_4)-3''$, \\
${\rm{sp}} \langle 3,2,1 \rangle ({\rm{T}}_5)-2'$,  ${\rm{sp}} \langle 3,2,5 \rangle ({\rm{T}}_5)-2'$, ${\rm{sp}} \langle 3,2,5 \rangle ({\rm{T}}_5)-2''$, \\
${\rm{sp}} \langle 3,7,2 \rangle ({\rm{T}}_5)-7'$,  ${\rm{sp}} \langle 3,7,2 \rangle ({\rm{T}}_5)-7''$, ${\rm{sp}} \langle 8,7,6 \rangle ({\rm{T}}_5)-7'$, \\
${\rm{sp}} \langle 8,7,6 \rangle ({\rm{T}}_5)-7''$.
\\

However, fourteen of these are isomorphic in pairs as follows: \\

\noindent  ${\rm{sp}} \langle 1,7,6 \rangle ({\rm{T}}_1)-7'  \cong  {\rm{sp}}\langle  1,7,6 \rangle ({\rm{T}}_1)-7''$,  
\\ ${\rm{sp}} \langle 6,8,1 \rangle ({\rm{T}}_2)-8'  \cong  {\rm{sp}}\langle  6,8,1 \rangle ({\rm{T}}_2)-8''$,  \\
 ${\rm{sp}} \langle 2,7,4 \rangle ({\rm{T}}_3)-7' \cong  {\rm{sp}}\langle  2,7,4 \rangle ({\rm{T}}_3)-7''$, \\ 
${\rm{sp}} \langle 7,8,4 \rangle ({\rm{T}}_3)-8' \cong  {\rm{sp}}\langle  3,2,1 \rangle ({\rm{T}}_5)-2'$, \\
 ${\rm{sp}} \langle 5,8,6 \rangle ({\rm{T}}_3)-8' \cong {\rm{sp}}\langle  5,8,6 \rangle ({\rm{T}}_3)-8''$, \\
${\rm{sp}} \langle 3,7,2 \rangle ({\rm{T}}_5)-7' \cong  {\rm{sp}}\langle  3,7,2 \rangle ({\rm{T}}_5)-7''$,    \\
  ${\rm{sp}} \langle 8,7,6 \rangle ({\rm{T}}_5)-7' \cong  {\rm{sp}}\langle  8,7,6 \rangle ({\rm{T}}_5)-7''$. \\

    \begin{center}
     \hspace*{0cm}      \pdfimage width 1.\textwidth {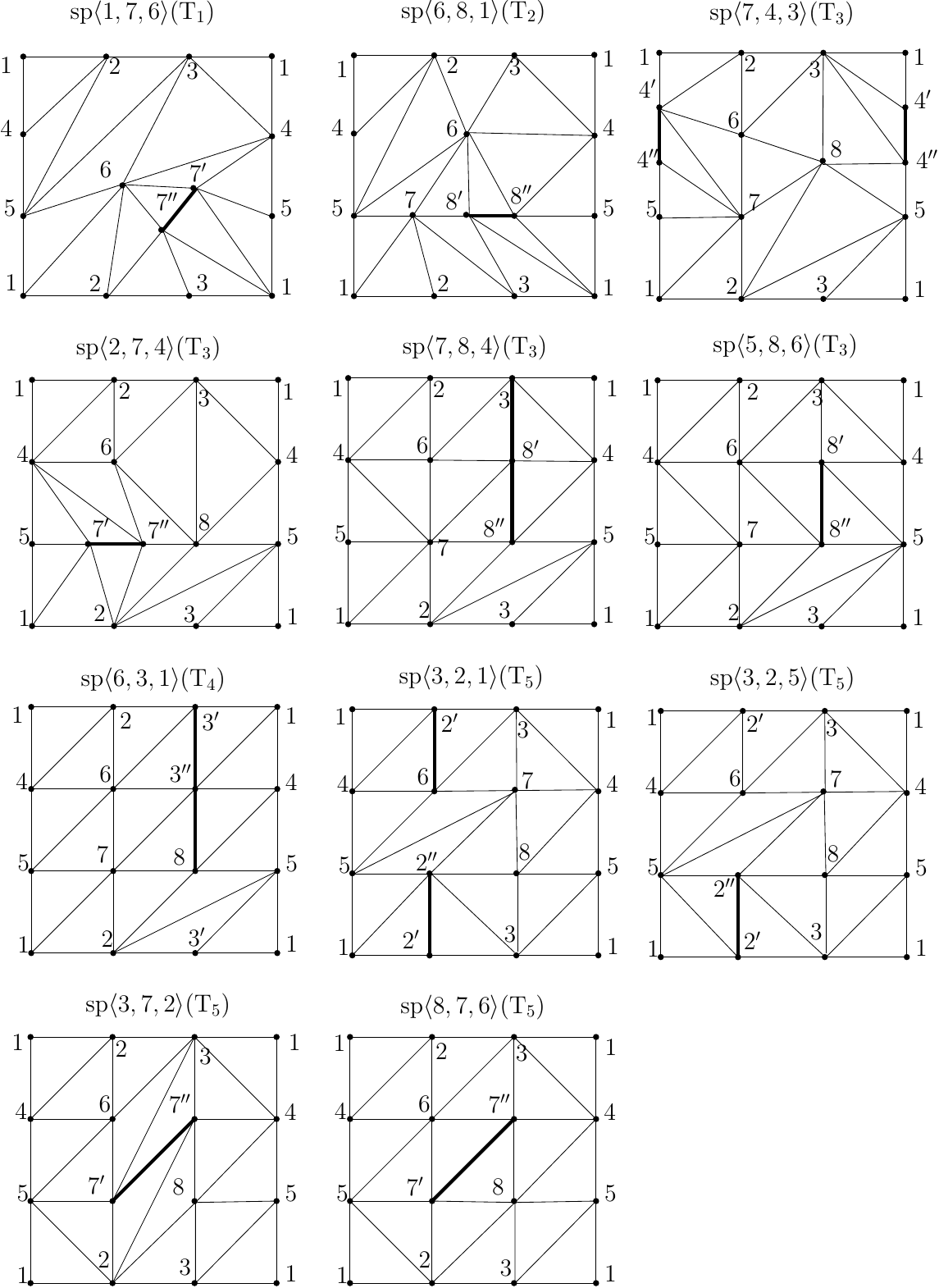}
\newline

    Figure 5. Pylonic triangulations in $\Lambda_3$.
    \end{center}

\medskip

We thus get Series 4 amounting to  $19-14/2=12$ non-isomorphic ITPTs:

 \begin{lemma}\label{10.2}
 There are precisely 12 ITPTs obtainable by deleting a  $\ast$-vertex from a  $\ast$-triangulation in  $\Lambda_3$.
 The twelve triangulations form Series 4 and are collected in Table~VI.
 \end{lemma}

\begin{proof}
Actually, it remains to verify that these twelve triangulations are non-isomor-
phic. This follows from
the fact that they have either differing d-vectors or bd-sequences; check with Table~VI.
\end{proof}

\medskip
    \begin{center}
 Table VI. ITPTs: Series 4.
\newline
\newline
        \begin{tabular}{|c|c |c|c | }
            \hline No. &  Triangulation        & d-vector                        &  bd-sequence
            \\
            \hline
           1     & ${\rm{sp}} \langle 1,7,6 \rangle (\rm{T}_1)-7'$ & (0,1,2,4)   &  (4,6,5,5,6)           \\
            \hline
            2     & ${\rm{sp}} \langle 6,8,1 \rangle (\rm{T}_2)-8'$ & (0,1,2,5,0) &  (4,6,5,5,6) \\
            \hline
            3    & ${\rm{sp}} \langle 7,4,3 \rangle (\rm{T}_3)-4'$ & (0,3,0,5,0) &  (4,6,4,6,4,6)          \\
            \hline
            4     & ${\rm{sp}} \langle 7,4,3 \rangle (\rm{T}_3)-4''$ & (0,0,5,2,1) &  (5,6,5,5,6)            \\
            \hline
           5  & ${\rm{sp}} \langle 2,7,4 \rangle (\rm{T}_3)-7'$ & (0,2,2,2,2) &  (4,7,5,4,7)
            \\
            \hline
            6   & ${\rm{sp}} \langle 7,8,4 \rangle (\rm{T}_3)-8'$ & (0,2,2,2,2) &  (4,6,4,5,7)  
            \\
            \hline
            7    & ${\rm{sp}} \langle 5,8,6 \rangle (\rm{T}_3)-8'$ & (0,1,3,3,1) & (4,5,5,6,6)           \\
            \hline
            8    & ${\rm{sp}} \langle 6,3,1 \rangle (\rm{T}_4)-3''$ & (0,2,4,0,2) & (4,5,5,4,5,5)         \\
            \hline
            9   & ${\rm{sp}} \langle 3,2,5 \rangle (\rm{T}_5)-2' $  & (1,2,1,2,2)&  (3,6,4,7,4,7)        \\
            \hline
            10   & ${\rm{sp}} \langle 3,2,5 \rangle (\rm{T}_5)-2''$ & (0,2,3,0,3)&   (4,5,7,5,7)       \\
            \hline
            11   & ${\rm{sp}} \langle 3,7,2 \rangle (\rm{T}_5)-7'$ & (0,3,1,1,3)& (4,6,7,4,7)
            \\
            \hline
            12  & ${\rm{sp}} \langle 8,7,6 \rangle (\rm{T}_5)-7'$ & (0,2,2,2,2) & (4,5,6,6,5)   \\
            \hline
        \end{tabular}
           \end{center}
\medskip

\medskip

  It is possible to produce more  $\ast$-triangulations by further splitting the  $\ast$-triangula-
tions in Lemma \ref{10.1}
  that have a unique rope, such as   ${\rm{sp}} \langle 6,8,1 \rangle ({\rm{T}}_2)$ (Figure 5). Some of the twice-split triangulations
   are out of  $\Lambda$, but some of them may have a rope whose contraction yields an IT; for instance,
$$ 
{\rm{sh}} \rangle 6,p\langle  \left(\begin{array}{c}
  {\rm{sp}}\langle  6,8',1 \rangle\left(\begin{array}{c}{\rm{sp}}\langle  6,8,1 \rangle ({\rm{T}}_2)\end{array}\right)\end{array}\right)  \cong {\rm{T}}_{16},
$$
\noindent where  $p$ stands for the  $\ast$-vertex.
    Thus the latter twice-split triangulations remain in  $\Lambda$!

 \begin{lemma}\label{10.3}
Splitting \ a \ $\ast$-triangulation \  $T^*$\  in \ Lemma \ref{10.1} \  (Figure 5) \ 
still produces a \ $\ast$-triangulation if and only if the splitting is equivalent to: (a) the
stellar subdivision of either of the two faces incident with a rope provided that rope is a unique
rope of  $T^*$, or (b) the cracking of a rope provided that rope is a unique rope of  $T^*$.
 \end{lemma}

\begin{proof}
The ``if'' part is obvious. In proving the ``only if'' part, observe from Figure 5 that there are two cases to consider as follows.

      {\it Case  }$``T^*$ {\it has only one rope''}. Observe from Figure 5 that the degrees
      of the end vertices of the rope are $5$ or $6$. Furthermore, splitting any divider centered at a $5$ or $6$-valent
       end vertex certainly destroys the pylonicity of   $T^*$ unless it meets condition (a) or (b); the only questionable
       situation is if the center is $6$-valent and the divider has spread $3$; then an additional consideration is as follows.
       There are precisely two $6$-valent end vertices of a rope---vertex  $4'$ in  ${\rm{sp}} \langle 7,4,3 \rangle ({\rm{T}}_3)$ and vertex $2'$  in
       ${\rm{sp}} \langle 3,2,5 \rangle ({\rm{T}}_5)$. By inspection, we verify that the following six triangulations produced by splitting a $3$-divider
       are not pylonic:  
$$
{\rm{sp}} \langle 1,4',7 \rangle \left(\begin{array}{c}{\rm{sp}}\langle  7,4,3 \rangle ({\rm{T}}_3)\end{array}\right),
$$
$$
{\rm{sp}} \langle 2,4',4''\rangle\left(\begin{array}{c}{\rm{sp}}\langle  7,4,3 \rangle ({\rm{T}}_3)\end{array}\right),      
$$
$$
{\rm{sp}} \langle 6,4',3 \rangle \left(\begin{array}{c}{\rm{sp}}\langle  7,4,3 \rangle ({\rm{T}}_3)\end{array}\right), 
$$
$$
{\rm{sp}} \langle 6,2',5 \rangle\left(\begin{array}{c}{\rm{sp}}\langle  3,2,5 \rangle ({\rm{T}}_5)\end{array}\right), 
$$
$$
{\rm{sp}} \langle 4,2',2'' \rangle \left(\begin{array}{c}{\rm{sp}}\langle  3,2,5 \rangle ({\rm{T}}_5)\end{array}\right),
$$
$$
{\rm{sp}} \langle 1,2',3 \rangle \left(\begin{array}{c}{\rm{sp}}\langle  3,2,5 \rangle ({\rm{T}}_5)\end{array}\right).
$$

        {\it      Case  } $``T^*$ {\it has precisely two ropes''.} In fact, there are only three such triangulations in
        Figure~5:  $T^* \in \{  {\rm{sp}}\langle  7,8,4 \rangle ({\rm{T}}_3),\,  {\rm{sp}}\langle  6,3,1 \rangle ({\rm{T}}_4),\,{\rm{sp}}\langle  3,2,1 \rangle ({\rm{T}}_5) \}$. Observe from
        Figure 5 that each $T^*$  contains precisely two ropes---which we denote by $[p,y]$  and
        $[p,z]$---and also observe that the degree of the central $\ast$-vertex   $p$ is $5$ or $6$, and that  $\mid y,p,z \mid \geq 2$.
        Assume that splitting \  ${\rm{sp}}\langle u,v,w \rangle $ \  of \  $T^*$  produces a $\ast$-triangulation. Then, necessarily,  $v=p$, since otherwise the
         newly produced edge  $[v',v'']$ would not be incident with the $\ast$-vertex  $p$. Furthermore, it is not hard to prove
         that for preserving the pylonicity property, it is necessary that $\langle u,v,w \rangle =\langle u,p,w \rangle $  is a $3$-divider that does not cross
          $\langle y,p,z \rangle $ at  $p$ and is edge disjoint from  $\langle y,p,z \rangle $. Since  $\mid y,p,z \mid \geq 2$, such a situation is theoretically
          possible but requires the degree of  $p$ to be at least $7$.
\end{proof}

By Corollary \ref{4.6}, there are not any other $\ast$-triangulations that belong to $\Lambda_3$  or can be obtained from a member of  $\Lambda_3$ by splitting. By Corollary \ref{4.2}, the corresponding ITPTs are obtained from Figure 5: (a) by the removal of either of the two faces incident with a single rope, and (b) by the removal of a single rope. In case (a), we have to inspect the sixteen triangulations obtained by the face removal from the eight single-roped triangulations in Lemma \ref{10.1}. These sixteen are naturally paired with each other in Table VII. The triangulations in each pair
$\sharp 1-6$ are isomorphic, which is verified straightforwardly. The ones in pair
$\sharp 7$ are not isomorphic because they have differing bd-sequences (check with Table VII). Moreover, as shown in the next paragraph, the ones in pair
$\sharp 8$ are not isomorphic even though they have the same $d$-vector and the same bd-sequence.

To show that the triangulations  
$$
{\rm{sp}} \langle 3,2,5 \rangle ({\rm{T}}_5)- [5,2',2'']
$$ 
\noindent and 
$$
{\rm{sp}} \langle 3,2,5 \rangle ({\rm{T}}_5)- [3,2',2'']
$$  
\noindent in pair
$\sharp 8$ of Table VII are not isomorphic, we observe on one hand that any such isomorphism would fix the single rope  $[2',2'']$, swapping the apices $5$ and $3$. On the other hand, the degrees of the neighboring vertices are ordered differently in the links of vertices $5$ and $3$, and hence vertices $5$ and $3$  are non-similar in ${\rm{sp}} \langle 3,2,5 \rangle ({\rm{T}}_5)$  and hence no such isomorphism is possible.

  There are no more isomorphic pairs in Table VII except the above-mentioned six pairs because the rest of the pairs have differing d-vectors or bd-sequences (except pair $\sharp
   8$). Therefore, Table VII provides Series 5 of 10 non-isomorphic ITPTs.
   
 In case (b) we obtain eight more ITPTs from the single-roped triangulations in Figure 5 by the rope removal. However, four of them are isomorphic to the ones originated from  $\Lambda_2$  and already present in Series 3 (Table II):  \\
\indent ${\rm{sp}} \langle 6,8,1 \rangle ({\rm{T}}_2)- [8',8''] \cong {\rm{T}}_{16} -i $, \\
\indent$ {\rm{sp}} \langle 7,4,3 \rangle ({\rm{T}}_3)- [4',4''] \cong {\rm{T}}_{17} -o$, \\ 
\indent$ {\rm{sp}} \langle 5,8,6 \rangle ({\rm{T}}_3)- [8',8''] \cong {\rm{T}}_{14} -g$, \\
\indent$ {\rm{sp}} \langle 8,7,6 \rangle ({\rm{T}}_5)- [7',7''] \cong {\rm{T}}_{12} -d $. \\ The remaining four form Series 6 of ITPTs, collected in Table~VIII; those four are pairwise non-isomorphic because they have differing $d$-vectors.

{\it Case 4.}  Parent triangulation  $T^*$ is in  $\Lambda_4$ or can be obtained from a member of
$\Lambda_4$ by a sequence of splittings. 
\medskip

 By Eq. (1), if  $T^*={\rm{sp}}\langle  u,v,w \rangle ({\rm{T}}_i) \in \Lambda_4$, then $\langle u,v,w \rangle $  is a $4$-divider. A comprehensive
 list of non-similar $4$-dividers can be generated by a sort of inclusion-exclusion technique like that
 introduced in Section \ref{S9}. We omit the details here since it is not very hard to determine such a list
 directly, using Table I. There are, in total, 42 non-similar $4$-dividers in the list. They are collected in
 Table IX in the form of the corresponding splittings of the corresponding ITs. The resulting split-up triangulations
  collectively form the whole set $\Lambda_4$  (with isomorphic duplications). It is a matter of a routine inspection
   to verify the following.

\newpage

    \begin{center}
 Table VII. ITPTs: Series 5.
\newline 

        \begin{tabular}{|c|c|c|c|c|}
            \hline No.  & No.  &                                    &               &             \\
            of  & of   & Triangulation                      & d-vector      & bd-sequence \\
            ITPT & pair &                                    &               &             \\
            \hline  1   &  1   & ${\rm{sp}} \langle 1,7,6 \rangle (\rm{T}_1)-[6,7',7'']\cong$  &(0,0,2,4,2)    & (5,5,7)     \\
            &      & $\cong {\rm{sp}}\langle  1,7,6 \rangle (\rm{T}_1)-[1,7',7'']$ &               &              \\
            \hline  2   &  2   & ${\rm{sp}} \langle 6,8,1 \rangle (\rm{T}_2)-[6,8',8'']\cong$  & (0,0,2,5,2,0) & (5,5,7)      \\
            &      & $\cong {\rm{sp}}\langle  6,8,1 \rangle (\rm{T}_2)-[1,8',8'']$ &               &              \\
            \hline  3   &  3   & ${\rm{sp}} \langle 7,4,3 \rangle (\rm{T}_3)-[7,4',4'']\cong$  & (0,0,3,3,3,0) & (5,6,7)      \\
            &      & $\cong {\rm{sp}}\langle  7,4,3 \rangle (\rm{T}_3)-[3,4',4'']$ &               &              \\
            \hline  4   &  4   & ${\rm{sp}} \langle 2,7,4 \rangle (\rm{T}_3)-[4,7',7'']\cong$  & (0,0,4,3,0,2) & (5,5,8)      \\
            &      & $\cong {\rm{sp}}\langle  2,7,4 \rangle (\rm{T}_3)-[2,7',7'']$ &               &              \\
            \hline  5   &  5   & ${\rm{sp}} \langle 3,7,2 \rangle (\rm{T}_5)-[3,7',7'']\cong$  & (0,1,4,0,2,2) & (5,5,8)      \\
            &      & $\cong {\rm{sp}}\langle  3,7,2 \rangle (\rm{T}_5)-[2,7',7'']$ &               &              \\
            \hline  6   &  6   & ${\rm{sp}} \langle 8,7,6 \rangle (\rm{T}_5)-[6,7',7'']\cong$  & (0,1,2,2,4,0) & (5,5,6)      \\
            &      & $\cong {\rm{sp}}\langle  8,7,6 \rangle (\rm{T}_5)-[8,7',7'']$ &  &     \\
            \hline  7   &  7   & $ {\rm{sp}}\langle  5,8,6 \rangle (\rm{T}_3)-[6,8',8'']$       & (0,0,3,3,3,0) & (5,5,6)      \\
            \hline  8   &  7   & $ {\rm{sp}}\langle  5,8,6 \rangle (\rm{T}_3)-[5,8',8'']$       & (0,0,3,3,3,0) & (5,5,7)    \\
            \hline  9   &  8   & ${\rm{sp}} \langle 3,2,5 \rangle (\rm{T}_5)-[5,2',2'']$        & (0,1,3,2,1,2) & (5,8,6)      \\
            \hline  10  &  8   &  ${\rm{sp}} \langle 3,2,5 \rangle (\rm{T}_5)-[3,2',2'']$        & (0,1,3,2,1,2) & (5,8,6)       \\
            \hline
        \end{tabular}
           \end{center}

\medskip
\medskip

    \begin{center}
 Table VIII. ITPTs: Series 6.
\newline

        \begin{tabular}{|c|c|c|c|}
            \hline No. & Triangulation             & d-vector        & bd-sequence  \\
            \hline 1   & ${\rm{sp}} \langle 1,7,6 \rangle (\rm{T}_1)-[7',7'']$ & (0,2,0,4,2)     & (7,4,7,4)    \\
            \hline 2   & ${\rm{sp}} \langle 2,7,4 \rangle (\rm{T}_3)-[7',7'']$ & (0,2,2,3,0,2)   & (4,8,4,8)    \\
            \hline 3   & ${\rm{sp}} \langle 3,2,5 \rangle (\rm{T}_5)-[2',2'']$ & (0,2,3,1,1,2)   & (4,8,5,8)    \\
            \hline 4   & ${\rm{sp}} \langle 3,7,2 \rangle (\rm{T}_5)-[7',7'']$ & (0,3,2,0,2,2)   & (4,8,4,8)    \\
            \hline
        \end{tabular}
           \end{center}

\medskip

\medskip

\section{The search for ITPTs: Case 4}
\label{S11}

 \begin{lemma}
 \label{11.1} None of the triangulations in $\Lambda_4$  are pylonic.
 \end{lemma}

 Thus, by Lemma \ref{11.1} and Corollaries \ref{4.2}, \ref{4.6}, there are no ITPTs which can be derived in Case 4.

\medskip
      
\noindent{\bf \textit{Proof of Lemma \ref{4.5}:}}
 It is not hard to prove that a $\Delta$ may occur only after the second consecutive splitting of an IT and only if
 the first applied splitting produces a single rope,  $\varepsilon$. We know from the above that a single rope may only
 occur as described above in the proof of Lemma \ref{10.3} (case $``T^*$ has only one rope''), in which event the degrees
 of the end vertices of  $\varepsilon$ are necessarily equal to $5$ or $6$, but not both equal to $6$.
 Then, it can be easily seen that the second splitting may lead to a $ \Delta$ only if it is the splitting of the 3-divider
 that contains $\varepsilon$  and is centered at the 6-valent end vertex of  $\varepsilon$. There are only two such second splittings
 (check with Figure 5),  ${\rm{sp}} \langle 2,4',4'' \rangle \left(\begin{array}{c}
 {\rm{sp}}\langle  7,4,3 \rangle ({\rm{T}}_3)\end{array}\right)$ \ and \   ${\rm{sp}} \langle 4,2',2'' \rangle \left(\begin{array}{c}{\rm{sp}}\langle  3,2,5 \rangle ({\rm{T}}_5)\end{array}\right)$,
 but neither has a $ \Delta$.

   \begin{center}
 Table IX. $\Lambda_4$ with isomorphic duplications.
\newline
\medskip
    \begin{tabular}{|cccc|c|}
        \hline                      &        Triangulations           &                    &                     &  Total  \\
        \hline  ${\rm{sp}} \langle 4,8,3 \rangle (\rm{T}_{8})$  & ${\rm{sp}} \langle 5,8,9 \rangle (\rm{T}_{8})$  &  ${\rm{sp}} \langle 1,8,7 \rangle (\rm{T}_{8})$     &                     &   3  \\
        \hline  ${\rm{sp}} \langle 1,4,5 \rangle (\rm{T}_{9})$  & ${\rm{sp}} \langle 3,4,8 \rangle (\rm{T}_{9})$  &  ${\rm{sp}} \langle 7,4,6 \rangle (\rm{T}_{9})$     & ${\rm{sp}} \langle 9,4,2 \rangle (\rm{T}_{9})$  &   4  \\
        \hline  ${\rm{sp}} \langle 1,4,5 \rangle (\rm{T}_{10})$ & ${\rm{sp}} \langle 2,4,9 \rangle (\rm{T}_{10})$ &  ${\rm{sp}} \langle 8,4,3 \rangle (\rm{T}_{10})$    &                     &   3   \\
        \hline  ${\rm{sp}} \langle 1,3,2 \rangle (\rm{T}_{14})$ & ${\rm{sp}} \langle 4,3,8 \rangle (\rm{T}_{14})$ &  ${\rm{sp}} \langle 6,3,5 \rangle (\rm{T}_{14})$    &                     &   3 \\
        \hline  ${\rm{sp}} \langle 2,3,1 \rangle (\rm{T}_{15})$ & ${\rm{sp}} \langle 8,3,4 \rangle (\rm{T}_{15})$ &  ${\rm{sp}} \langle 9,3,7 \rangle (\rm{T}_{15})$    & ${\rm{sp}} \langle 5,3,6 \rangle (\rm{T}_{15})$ &   4  \\
        \hline  ${\rm{sp}} \langle 2,3,1 \rangle (\rm{T}_{16})$ & ${\rm{sp}} \langle 8,3,4 \rangle (\rm{T}_{16})$ &  ${\rm{sp}} \langle 5,3,6 \rangle (\rm{T}_{16})$    &                     &   3  \\
        \hline  ${\rm{sp}} \langle 1,3,2 \rangle (\rm{T}_{17})$ & ${\rm{sp}} \langle 4,3,8 \rangle (\rm{T}_{17})$ &  ${\rm{sp}} \langle 7,3,9 \rangle (\rm{T}_{17})$    & ${\rm{sp}} \langle 6,3,5 \rangle (\rm{T}_{17})$ &   4  \\
        \hline ${\rm{sp}} \langle 1,4,9 \rangle (\rm{T}_{18})$  & ${\rm{sp}} \langle 2,4,8 \rangle (\rm{T}_{18})$ & ${\rm{sp}} \langle 6,4,7 \rangle (\rm{T}_{18})$     & ${\rm{sp}} \langle 5,4,3 \rangle (\rm{T}_{18})$ &   4 \\
        \hline ${\rm{sp}} \langle 1,4,9 \rangle (\rm{T}_{19})$  & ${\rm{sp}} \langle 2,4,8 \rangle (\rm{T}_{19})$ & ${\rm{sp}} \langle 6,4,7 \rangle (\rm{T}_{19})$     & ${\rm{sp}} \langle 5,4,3 \rangle (\rm{T}_{19})$ &   4 \\
        \hline ${\rm{sp}} \langle 1,4,9 \rangle (\rm{T}_{20})$  & ${\rm{sp}} \langle 2,4,8 \rangle (\rm{T}_{20})$ & ${\rm{sp}} \langle 6,4,7 \rangle (\rm{T}_{20})$     & ${\rm{sp}} \langle 5,4,3 \rangle (\rm{T}_{20})$ &     \\
        ${\rm{sp}} \langle 9,3,7 \rangle (\rm{T}_{20})$  & ${\rm{sp}} \langle 5,3,6 \rangle (\rm{T}_{20})$ & ${\rm{sp}} \langle 1,3,2 \rangle (\rm{T}_{20})$     & ${\rm{sp}} \langle 4,3,8 \rangle (\rm{T}_{20})$ &   8  \\
        \hline ${\rm{sp}} \langle 2,10,3 \rangle (\rm{T}_{21})$ & ${\rm{sp}} \langle 6,10,8 \rangle (\rm{T}_{21})$&                         &                     &   2  \\
        \hline

    \end{tabular}
        \end{center}
\medskip

 \section{Concluding theorem}
 \label{S12}

Combining the results in the previous sections, we have identified a total of 297 non-isomorphic ITPTs as detailed in the following summarizing theorem.

\begin{theorem}
\label{12.1}
 Up to isomorphism, there are totally 297 irreducible triangulations of the once-punctured torus. They are presented in six series as follows: 80 triangulations in Series 1 (Section \ref{S6}), 129 triangulations in Series 2 (Section \ref{S7}), 62 triangulations in Series 3 (Section~\ref{S8}), and 12, 10, and 4 triangulations in (respectively) Series 4 in Table VI, Series 5 in Table VII, and Series~6 in Table VIII (all in Section~\ref{S10}).
 \end{theorem}

  \begin{proof}
  The considerations of Sections \ref{S6}--\ref{S11}, along with Corollaries \ref{4.2} and \ref{4.6}, guarantee that we have not missed any ITPT in the search. It remains to show that all the ITPTs we have found are pairwise non-isomorphic.

  We have verified in Section \ref{S6} that Series 1 contains no isomorphic pairs, nor
do Series 2 (Section \ref{S7}), Series 3 (Section \ref{S8}), nor any of Series 4--6 (Section \ref{S10}). It remains to prove that Series 1--6 are pairwise disjoint from each other.

By Corollary \ref{5.2}, \ Series 2 \  (obtained by deleting the $\ast$-vertex from the members of  $\Lambda_1$), Series 3 (from  $\Lambda_2$), and Series 4 (from  $\Lambda_3$) are pairwise disjoint. Clearly, each of these three Series is disjoint from Series 1 (produced immediately from the toroidal ITs). Furthermore, each of the triangulations in Series 5 and 6, with the patch restored, can be produced from some IT by two, but not by one, consecutive splittings, hence the triangulations resulting from these splittings collectively are out of  $\Lambda$, and hence Series 5 and 6 are disjoint from the other Series.
Finally, the triangulations in the union of Series 5 and 6 have pairwise differing d-vectors or bd-sequences (except the non-isomorphic pair $\sharp 8$ in Series 5 considered in Section \ref{S10}) and thus are all non-isomorphic.
   \end{proof}

{ \bf Acknowledgment.} The first author (Lawrencenko) is grateful to Professor Branko Gr\"{u}nbaum for enlightening correspondence on the 1934 original German edition of the book \cite{SR}.

\bibliographystyle{model1-num-names}
\bibliography{<your-bib-database>}

\newpage

\noindent \textsc{S. Lawrencenko} \\
\noindent Russian State University of Tourism and Service, \\
\noindent Institute of Service Technologies,\\
\noindent 20 Krasnaya, Podolsk, Moscow Region, 142116, Russia\\
\noindent e-mail: \url{lawrencenko@hotmail.com}\\

\noindent \textsc{T. Sulanke}\\
\noindent Department of Physics, Indiana University,\\
\noindent Bloomington, Indiana 47405, USA\\
\noindent e-mail:  \url{tsulanke@indiana.edu}\\

\noindent \textsc{M.\,T. Villar}\\
\noindent Dpto. de Geometr\'{\i}a y Topolog\'{\i}a, Universidad de Sevilla,\\
\noindent C/ Tarfia s/n, 41012 Sevilla, Spain\\
\noindent e-mail:  \url{villar@us.es}\\

\noindent \textsc{L.\,V. Zgonnik} \\
\noindent Russian State University of Tourism and Service, \\
\noindent Institute for Tourism and Hospitality,\\
\noindent 32A Kronstadt Boulevard, Moscow, 125438, Russia\\
\noindent e-mail:  \url{mila.zgonnik1@yandex.ru}\\

\noindent \textsc{M.\,J. Ch\'avez}\\
\noindent Dpto. de Matem\'atica Aplicada I, Universidad de Sevilla,\\
\noindent Avda. Reina Mercedes s/n, Sevilla, Spain\\
\noindent e-mail:  \url{mjchavez@us.es}\\

\noindent \textsc{J.\,R. Portillo}\\
\noindent Dpto. de Matem\'atica Aplicada I, Universidad de Sevilla,\\
\noindent Avda. Reina Mercedes s/n, Sevilla, Spain\\
\noindent e-mail: \url{josera@us.es}\\

\end{document}